\documentclass[11pt]{amsart}

\usepackage[active]{srcltx} %SRC Specials for DVI Searching
\usepackage{amsmath}
\usepackage{amsfonts}
\usepackage{amssymb}
\usepackage{latexsym}
\usepackage{graphicx}
\usepackage[cp866]{inputenc}

 \advance\hoffset by -0.0 truecm

\newtheorem{theorem}{Theorem}[section]
\newtheorem{corollary}[subsection]{Corollary}
\newtheorem{lemma}[subsection]{Lemma}
\newtheorem{proposition}[subsection]{Proposition}
\theoremstyle{definition}
\newtheorem{definition}[subsection]{Definition}
\theoremstyle{remark}
\newtheorem{remark}[theorem]{Remark}
\numberwithin{equation}{section}

\newcommand{\R}{\mathbb{R}}
\newcommand{\C}{\mathbb{C}}
\newcommand{\ii}{\mathbf{i}}
\newcommand{\jj}{\mathbf{j}}
\newcommand{\kk}{\mathbf{k}}
\newcommand{\I}{\mathcal{M}_1}
\newcommand{\J}{\mathcal{M}_2}
\newcommand{\K}{\mathcal{M}_3}
\newcommand{\Ib}{\mathbf{M_1}}
\newcommand{\Jb}{\mathbf{M_2}}
\newcommand{\Kb}{\mathbf{M_3}}
\newcommand{\M}{\mathcal{M}}
\newcommand{\Mb}{\mathbf{M}}
\newcommand{\nb}{\mathbf{n}}

\newcommand{\Tb}{\mathbf{\Theta}}
\newcommand{\bx}{\mathbf{x}}
\newcommand{\by}{\mathbf{y}}
\newcommand{\T}{\mathcal{T}}
\newcommand{\U}{\mathcal{U}}

\newcommand{\re}{\hbox{Re}}
\newcommand{\im}{\hbox{Im}}

\begin{document}

\title[Quaternion $\mathbb H$-type groups $Q^n$]
{ Anisotropic quaternion Carnot groups: geometric analysis and
Green's function}

\author{Der-Chen Chang,\ Irina Markina}

\address{Department of Mathematics, Georgetown University, Washington
D.C. 20057, USA}

\email{chang@georgetown.edu}

\address{Department of Mathematics,
University of Bergen, Johannes Brunsgate 12, Bergen 5008, Norway}

\email{irina.markina@uib.no}
\thanks{This work was supported by
Projects FONDECYT (Chile) \# 7050181, \#1040333, and by grant of
the University of Bergen}

\subjclass[2000]{53C17, 53C22, 35H20}

\keywords{Quaternions, Carnot-Carath\'eodory metric, nilpotent Lie
groups, Hamiltonian formalism, Green function.}

%%% ----------------------------------------------------------------------

\begin{abstract}
We construct examples of $2$-step Carnot groups related to
quaternions and study their fine structure and geometric
properties. This involves the Hamiltonian formalism, which is used
to obtain explicit equations for geodesics and the computation of
the number of geodesics joining two different points on these
groups. We able to find the explicit lengths of geodesics. We
present the fundamental solutions of the Heat and sub-Laplace
equations for these anisotropic groups and obtain some estimates
for them, which may be useful.
\end{abstract}

%%% ----------------------------------------------------------------------
\maketitle
%%% ----------------------------------------------------------------------

\section{Introduction}

This paper presents examples of Carnot groups and studies their
fine structure, geometric properties and basic differential
operators attached to them. A Carnot group is a connected and
simply connected $m$-step nilpotent Lie group $\mathbb G$ whose
Lie algebra $\mathcal G$ decomposes into the direct sum of vector
subspaces $V_1\oplus V_2\oplus\ldots\oplus V_m$ satisfying the
following relations:
$$[V_1,V_k] = V_{k+1},\qquad 1\leq k<m,\quad\quad [V_1,V_m]
=\{0\}.$$ The simplest examples of this group are Euclidean space
$\mathbb R^n$, Heisenberg group $\mathbb H^n$ and $\mathbb
H$(eisenberg)-type groups introduced by Kaplan~\cite{Kap1}. The
Carnot groups form a natural habitat for extensions of many of the
objects studied in Euclidean space and find applications in the
study of strongly pseudoconvex domains in complex analysis,
semiclassical analysis of quantum mechanics, control theory,
probability theory of degenerate diffusion process and others. The
geometry of Carnot groups and differential operators related with
them were studied extensively by many mathematicians, for
instance, in~\cite{BMT,BGGR3,CChGr3,CG,FS,Ga,Kap2,Kor1,MShS}.

We construct examples of $2$-step Carnot groups, related to the
multidimensional space of quaternion numbers. We will call these
groups {\it anisotropic quaternion groups} and denote them by~
$Q^n$. In~\cite{ChM} the quaternion $\mathbb H$-type groups were
studied. The results of~\cite{ChM} can be easily extended to the
multidimensional quaternion space. The examples of the present
paper contain the multidimensional quaternion $\mathbb H$-type
groups as a particular case. We construct the Hamiltonian function
associated with the sub-Laplacian generated by left invariant
vector fields. Solving the Hamiltonian system of differential
equations we give exact solutions that describe the geodesics on
the group. We study geodesic connectivity between any two points
of the group. It is known that every point of a Riemannian
manifold is connected to every other point in a sufficiently small
neighborhood by one single, unique geodesic. But in this case,
there will be points arbitrarily near a point which are connected
to this point by an infinite number of geodesics. Since we are
working on a group, we may simply assume that the point is the
origin $O=(0,0)$. We prove the following results:

{\it $(1)$. If $P=(x,0)=(x_1,\ldots, x_n,0)$ with $x\ne 0$, then
there is only one geodesic connecting the origin $O$ and the point
$P$;

$(2)$. If $P=(x_1,\ldots, x_n,z)$ with $x_l\ne 0$, $l=1,\ldots,n$
and $z\ne 0$, then there are finitely many geodesics connecting
the origin $O$ and the point $P$;

$(3)$. If $P=(x_1,\ldots,x_n,z)$ with $x_l\ne 0$ for $l=1,\ldots,
p-1$ and $x_l=0$ for $l=p,\ldots, n$, $z\ne 0$, then there are
countably infinitely many geodesics connecting the point $O$ and
the point $P$;

$(4)$. If $P=(0,z)$ with $z\ne 0$, then there are uncountably
infinitely many geodesics connecting the point $O$ and the point
$P$.}

We will discuss basic properties of geodesics in sections 2 to 4.
Then we will prove connectivity theorems in section 5 (see
Theorems \ref{54}, \ref{4}, \ref{3}, and \ref{an19}). Furthermore,
parametric equations and arc lengths of all these geodesics are
calculated explicitly.

We also consider complex geodesics and find a relation between the
complex action function and the Carnot-Carath\'{e}odory metric.
The complex action function allows us to deduce the transport
equation and its solution: volume element. The fundamental
solution of the Heat equation is given in terms of the complex
action function and volume element. More precisely, the heat
kernel at the origin is given by
$$P(y,w,t)=\frac{C}{t^{2n+3}}\int_{\mathbb
R^3}e^{\frac{-f}{t}}V(\tau)\,d\tau,$$ where $$f(y,w,\tau) =
-i\sum_m\tau_mw_{m}+\sum_{l=1}^{n}\frac{|y_l|^2}{4}|\tau|_l\coth(|\tau|_l)$$
is the modified complex action and
$$V(\tau)=\prod_{l=1}^{n}\frac{|\tau|_l^{2}}{\sinh^{2}(|\tau|_l)}$$ is the volume
element. (See Theorem \ref{th:heat1}). Integrating the fundamental
solution of the Heat equation with respect to the time variable,
we obtain the Green function for the sub-Laplacian (see Theorem
\ref{th:fund5}):
\[
G(x,z)=-\frac{2^{2n}(2\pi)^{2n+3}}{(2n+1)!}\int_{{\mathbb
R}^3}\frac{V(\tau+i\varepsilon
\tilde{z})}{f^{2n+2}(\tau+i\varepsilon \tilde{z})}d\tau.
\]
The last section is devoted to some estimates of fundamental
solutions, which may be useful.

Part of this paper was finished while the authors visited
Universidad T\'ecnica Federico Santa Mar\'{i}a, Valpara\'{\i}so,
Chile in December, 2005 under the grant Projects FONDECYT  (Chile)
\# 7050181, \#1040333. We would like to thank Professor Alexander
Vasil'ev and the Departamento de Matem\'atica of UTFSM for their
invitation and the warm hospitality extended to them during their
stay in Chile. We would also like to thank Professor Eric Grinberg
for many inspired conversations on this project.

\section{Definitions}

A quaternion is a mathematical concept (re)introduced by William
Rowan Hamilton from Ireland in 1843 \cite{B}. (It has been said that
{\it when Hamilton discovered the quaternions, they stayed
discovered}).  The idea captured the popular imagination for a time
because it involved relatively simple calculations that abandon the
commutative law, one of the basic rules of arithmetic. Specifically,
a quaternion is a non-commutative extension of the complex numbers.
As a vector space over the real numbers, the quaternions have
dimension $4$, whereas the complex numbers have dimension $2$. While
the complex numbers are obtained by adding the element $\ii$ to the
real numbers which satisfies $\ii^2=-1$, the quaternions are
obtained by adding the elements $\ii,\jj$, and $\kk$ to the real
numbers which satisfy the following relations
$$\ii^2=\jj^2=\kk^2=\ii\jj\kk=-1.$$
Unlike real or complex numbers, multiplication of quaternions is
not commutative, e.~g.,\begin{equation}\label{12} \ii\jj=-
\jj\ii=\kk,\quad \jj\kk=- \kk\jj=\ii,\quad \kk\ii=- \ii\kk=\jj.
\end{equation} The quaternions are an example of a division ring,
an algebraic structure similar to a field except for commutativity
of multiplication. In particular, multiplication is still
associative and every non-zero element has a unique inverse.

The quaternions can be written  as a combination of a scalar and a
vector in analogy with the complex numbers being representable as a
sum of real and imaginary parts, $a\cdot 1+b\cdot\ii$. For a
quaternion $h=a+b\ii+c\jj+d\kk$ we call a scalar $a$  the {\it real
part} and the $3$-dimensional vector $\mathbf u=b\ii+c\jj+d\kk$ is
called the {\it imaginary part} of $h$ and is a {\it pure
quaternion}. In $\mathbb R^4$, the basis of quaternion numbers can
be given by real matrices
\begin{equation*}\mathcal U=\left[\array{rrrr}
1 & 0 & 0 & 0
\\
0 & 1 & 0 & 0
\\
0 & 0 & 1 & 0
\\
0 & 0 & 0 & 1\endarray\right],\end{equation*}\begin{equation*}
\I=\left[\array{rrrr}0 & 1 & 0 & 0
\\ -1 & 0 & 0 & 0
\\
0 & 0 & 0 & 1
\\ 0 & 0 & -1 & 0\endarray\right],\qquad
\J=\left[\array{rrrr} 0 & 0 & 0 & -1
\\
0 & 0 & -1 & 0
\\
0 & 1 & 0 & 0
\\
1 & 0 & 0 & 0\endarray\right], \qquad \K=\left[\array{rrrr} 0 & 0
& -1 & 0
\\
0 & 0 & 0 & 1
\\
1 & 0 & 0 & 0
\\
0 & -1 & 0 & 0\endarray\right].\end{equation*} We have
$$h=\left[\array{rrrr}
a & b & -d & -c
\\
-b & a & -c & d
\\
d & c & a & b
\\
c & -d & -b & a\endarray\right]=a\mathcal U+b\I+c\J+d\K.$$ Similarly
to complex numbers, vectors, and matrices, the addition of two
quaternions is equivalent to summing up the coefficients. Set
$h=a+\mathbf u$, and $q=t+x\ii+y\jj+z\kk=t+\mathbf v$. Then
$$h+q=(a+t)+(\mathbf u+\mathbf v)=(a+t)+(b+x)\ii+(c+y)\jj+(d+z)\kk.$$
Addition satisfies all the commutation and association rules of real
and complex numbers. The quaternion multiplication (the Grassmanian
product) is defined by $$hq=(at-\mathbf u\cdot\mathbf v)+(a\mathbf
v+t\mathbf u+\mathbf u\times\mathbf v),$$ where $\mathbf
u\cdot\mathbf v$ is the scalar product and $\mathbf u\times\mathbf
v$ is the vector product of $\mathbf u$ and $\mathbf v$, both in
$\mathbb R^3$. The multiplication is not commutative because of the
non-commutative vector product. The non-commutativity of
multiplication has some unexpected consequences, e.g., polynomial
equations over the quaternions may have more distinct solutions then
the degree of a polynomial. The equation $h^2+1=0$, for instance,
has infinitely many quaternion solutions $h=a+b\ii+c\jj+d\kk$ with
$b^2+c^2+d^2=1$. The {\it conjugate} of a quaternion
$h=a+b\ii+c\jj+d\kk$, is defined as $h^*=a-b\ii-c\jj-d\kk$ and the
{\it absolute value} of $h$ is defined as $|h|=\sqrt{h
h^{*}}=\sqrt{a^2+b^2+c^2+d^2}$.

Let us denote the space of quaternions by $\mathcal H$. We
consider $n$-tuples of quaternions: $H=(h_1,h_2,\ldots,h_n)$,
$h_i\in\mathcal H$, $i=1,\ldots,n$. We may define addition between
two of them in an obvious way:
$$H+Q=(h_1+q_1,h_2+q_2,\ldots,h_n+q_n)$$
for $H=(h_1,\ldots, h_n)$ and $Q=(q_1,\ldots,q_n)$. Multiplication
by scalar is defined as $\alpha H=(\alpha h_1,\alpha
h_2,\ldots,\alpha h_n)$, where $\alpha$ may be real or complex
number. Therefore, the $n$-dimensional quaternions,  $\mathcal
H^n$ is a vector space. The norm is defined by
$$|H|_{\mathcal H^n}=\Big(\sum_{i=1}^{n}|h_i|^2\Big)^{1/2}.$$

In this article we will construct $2$-step Carnot groups related to
the multidimensional quaternion numbers. We will call these groups
{\it anisotropic quaternion groups} $Q^n$. To explain precisely the
idea of their construction we first introduce the notion of $\mathbb
H$-type Carnot groups and then, making some modifications, we arrive
at our main example.

$\mathbb H$-type homogeneous groups are simply connected $2$-step
Lie groups $\mathbb G$ whose algebras $\mathcal G$ are graded and
carry an inner product such that
\begin{itemize}
\item[(i)]{$\mathcal G$ is the orthogonal direct sum of the
generating subspace $V_1$ and the center $V_2$: $$\mathcal
G=V_1\oplus V_2,\quad \ V_2=[V_1,V_1],\quad  [V_1,V_2]=0,$$}
\item[(ii)]{the homomorphisms $J_{\mathcal Z}:V_1\to V_1$,
${\mathcal Z}\in V_2$ defined by
$$\langle J_{\mathcal Z}X,X^{\prime}\rangle=\langle {\mathcal Z},[X,X^{\prime}]\rangle,
\quad X,X^{\prime}\in V_1$$ satisfy the equation
$$J^2_{\mathcal Z}=-|\mathcal Z|^2I, \quad {\mathcal Z}\in V_2.$$}
\end{itemize}
Here $\langle\cdot,\cdot\rangle$ is a positive definite
non-degenerating quadratic form on $\mathcal G$, $[\cdot,\cdot]$ is
a commutator and $I$ is the identity. The group is generated from
its algebra by exponentiation.

To construct the multidimensional quaternion $\mathbb H$-type group
we take the space of quaternions $\mathcal H^n$ as $V_1$ and
generate the center $V_2$. We consider the $n$-dimensional imaginary
quaternions $${\mathcal Z}_1=(a_1\ii,\ldots,a_1\ii),\ \ {\mathcal
Z}_2=(a_2\jj,\ldots,a_2\jj),\ \ {\mathcal
Z}_3=(a_3\kk,\ldots,a_3\kk)$$ with positive constants $a_m$. They
have the following representation as real matrices $4n\times 4n$:
\begin{equation*}
M_m=\left[\array{ccc}a_m\M_m & \ & 0
\\ \ & \ddots & \
\\
 0  & \ &a_m\M_m
\endarray\right],\end{equation*}
where there are $n$ blocks on the
diagonal of each matrix $M_m$, $m=1,2,3$. The matrices $M_m$,
$m=1,2,3$, are the matrices associated to the homomorphisms
$J_{\mathcal Z}$.

Now we extend the construction, introducing anisotropy to this very
symmetric setting. We take an arbitrary $n$-dimensional imaginary
quaternions $${\mathcal Z}_1=(a_{11}\ii,\ldots,a_{1n}\ii),\ \
{\mathcal Z}_2=(a_{21}\jj,\ldots,a_{2n}\jj),\ \ {\mathcal
Z}_3=(a_{31}\kk,\ldots,a_{3n}\kk),$$ with $a_{ml}>0$ for all
$m=1,2,3$ and $l=1,\ldots,n$. The representation as real matrices
$4n\times 4n$ are following:
\begin{equation*}
\Ib=\left[\array{ccc}a_{11}\I & \ & 0
\\ \ & \ddots & \
\\
 0  & \ &a_{1n}\I
\endarray\right],\qquad
\Jb=\left[\array{ccc} a_{21}\J & \ & 0
\\
\ & \ddots & \
\\
0 & \  & a_{2n}\J
\endarray\right], \end{equation*}
\begin{equation*}
\Kb=\left[\array{ccc} a_{31}\K & \ & 0
\\
\  & \ddots & \ \\
 0 & \  & a_{3n}\K
\endarray\right],\end{equation*} where there are $n$ blocks on the
diagonal of each matrix $\Mb_m$, $m=1,2,3$. We construct and make
the principal calculations for the anisotropic quaternion group
$Q^n$ with center $V_2$ of topological dimension $3$. For the other
examples see Remark~\ref{r1}. The corresponding algebra ${\mathcal
Q}^n$ is the two-step algebra $V_1\oplus V_2$. The topological
dimensions of the group is $\dim Q^n=4n+3$. The homogeneous
dimension defined by the formula $\nu=\dim V_1+2\dim V_2$ plays an
important role in analysis on homogeneous groups. We see that the
homogeneous dimension is always greater than the topological
dimension and in our case equals $\nu(Q^n)=4n+6$. An isotropic case
based on the one dimensional space of quaternions was studied
in~\cite{ChM} and the anisotropic Heisenberg group was considered
in~\cite{CChGr3}.

We set the standard orthonormal systems $\{X_{kl}\}\in \mathbb
R^{4n}$, $k=1,2,3,4$, $l=1,\ldots,n$ and $\{Z_{1},Z_{2},Z_{3}\}\in
\mathbb R^3$. We reserve the following indexes: $l=1,\ldots n$
denotes the coordinate index in $(h_1,\ldots,h_n$); $k=1,2,3,4$
denotes the index of coordinates inside of each quaternion $h_l$ and
$m=1,2,3$ is related to the coordinate index in $V_2$ or the index
of the matrices $\Mb_m$. The matrices $\Mb_m$ transform the basis
vectors in the following way
\begin{equation}\label{63}\begin{array}{llll} \Mb_1 X_{1l}  = -a_{1l}X_{2l}, \quad
& \Mb_1 X_{2l}  =  a_{1l}X_{1l},\quad  & \Mb_1 X_{3l}  =
-a_{1l}X_{4l},\quad & \Mb_1 X_{4l} =  a_{1l}X_{3l},
\\
\Mb_2 X_{1l}  =  a_{2l}X_{4l},\quad & \Mb_2 X_{2l}  =
a_{2l}X_{3l},\quad & \Mb_2 X_{3l} =  -a_{2l}X_{2l},\quad & \Mb_2
X_{4l} = -a_{2l}X_{1l},
\\
\Mb_3 X_{1l}  =  a_{3l}X_{3l},\quad & \Mb_3 X_{2l}  =
-a_{3l}X_{4l},\quad & \Mb_3 X_{3l}  =  -a_{3l}X_{1l},\quad & \Mb_3
X_{4l} =  a_{3l}X_{2l}.\end{array}\end{equation}

We use the normal coordinates
$$q=(x,z)=(x_{11},x_{21},x_{31},x_{41},\ldots,x_{1n},x_{2n},x_{3n},x_{4n},z_{1},z_{2},z_{3})$$ for the elements
$$\exp\Big(\sum_{k,l}x_{kl}X_{kl}
+\sum_{m}z_{m}Z_{m}\Big)\in Q^n.$$ The Baker-Campbell-Hausdorff
formula
$$\exp(X+Z)\exp(X^{\prime}+Z^{\prime})=
\exp(X+X^{\prime},Z+Z^{\prime}+\frac{1}{2}[X,X^{\prime}]),$$ for
$X,X^{\prime}\in V_1$, $Z,Z^{\prime}\in V_2$ defines the
multiplication law on $Q^n$. Precisely,  we have
\begin{eqnarray*}
L_q(q^{\prime}) & = & L_{(x,z)}(x^{\prime},z^{\prime}) =
(x,z)\circ(x^{\prime},z^{\prime}) \\ & = & \big(x+x^{\prime},
z_{1}+z_{1}^{\prime}+\frac{1}{2}(\mathbf M_1 x, x^{\prime}),
z_{2}+z_{2}^{\prime}+\frac{1}{2}(\mathbf M_2 x, x^{\prime}),
z_{3}+z_{3}^{\prime}+\frac{1}{2}(\mathbf M_3 x, x^{\prime})\big),
\end{eqnarray*} for
$q=(x,z)$ and $q^{\prime}=(x^{\prime},z^{\prime})$, where
$(\mathbf M_m x,x^{\prime})$ is the usual scalar product of the
vector $\mathbf M_m x\in\mathbb R^{4n}$ by $x^{\prime}\in\mathbb
R^{4n}$. The multiplication ``$\circ$" defines the left
translation $L_q$ of $q^{\prime}=(x^{\prime},z^{\prime})$ by the
element $q=(x,z)$ on the group $Q^n$.

We associate the Lie algebra ${\mathcal Q}^n$ of the group $Q^n$
with the set of all left invariant vector fields of the tangent
bundle $TQ^n$. The tangent bundle contains a natural subbundle
${\T}Q^n$ consisting of ``horizontal" vectors. We call $\T Q^n$
the {\it horizontal bundle}. The horizontal bundle is spanned by
the left-invariant vector fields $\widetilde
X_{11}(x,z),\ldots,\widetilde X_{4n}(x,z)$ with $\widetilde
X_{kl}(0,0)=X_{kl}$, $k=1,\ldots,4$, $l=1,\ldots,n$, (see for
example,~\cite{CChGr1}, \cite{CChGr3}, ~\cite{Reim}). In
coordinates of the standard Euclidean basis
$\frac{\partial}{\partial x_{kl}}$, $\frac{\partial}{\partial
z_{m}}$, these vector fields are expressed as
\begin{equation}
\label{62}
\begin{split} \widetilde X_{1l}(x,z)  = & \frac{\partial}{\partial
x_{1l}}+\frac{1}{2}\Big(+a_{1l}x_{2l}\frac{\partial}{\partial
z_{1}}-a_{2l}x_{4l}\frac{\partial}{\partial z_{2}}-
a_{3l}x_{3l}\frac{\partial}{\partial z_{3}}\Big), \\
\widetilde X_{2l}(x,z) = & \frac{\partial}{\partial
x_{2l}}+\frac{1}{2}\Big(-a_{1l}x_{1l}\frac{\partial}{\partial
z_{1}}-a_{2l}x_{3l}\frac{\partial}{\partial z_{2}}+
a_{3l}x_{4l}\frac{\partial}{\partial z_{3}}\Big),\\
\widetilde X_{3l}(x,z) = & \frac{\partial}{\partial
x_{3l}}+\frac{1}{2}\Big(+a_{1l}x_{4l}\frac{\partial}{\partial
z_{1}}+a_{2l}x_{2l}\frac{\partial}{\partial z_{2}}+
a_{3l}x_{1l}\frac{\partial}{\partial z_{3}}\Big), \\
\widetilde X_{4l}(x,z) = & \frac{\partial}{\partial
x_{4l}}+\frac{1}{2}\Big(-a_{1l}x_{3l}\frac{\partial}{\partial
z_{1}}+a_{2l}x_{1l}\frac{\partial}{\partial z_{2}}-
a_{3l}x_{2l}\frac{\partial}{\partial z_{3}}\Big),
\end{split}
\end{equation} for $l=1,\ldots,n$. The left invariant vector fields $\widetilde
Z_{m}(x,z)$ with $\widetilde Z_{m}(0,0)=Z_{m}$, $m=1,2,3$, are
simply the vector fields \begin{equation} \label{62z}\widetilde
Z_{m}(x,z)=\frac{\partial}{\partial z_{m}}.\end{equation} We write
simply $X_{kl}$ and $Z_{m}$ instead of $\widetilde X_{kl}(x,z)$
and $\widetilde Z_{m}(x,z)$, if no confusion may arise. Note that
if we fix $m$ to be equal only to $1,2$ or $3$ then the vector
fields~\eqref{62} are reduced to the anisotropic vector fields of
anisotropic $\mathbb H^{2n}$ Heisenberg group and the group $Q^n$
is isomorphic to anisotropic $\mathbb H^{2n}$ Heisenberg group,
considered in~\cite{CChGr3}. We also use the notation
$X_l=(X_{1l},\ldots, X_{4l})$, $l=1,\ldots,n$. We call the next
vector $$X  = (X_{11},\ldots,X_{4\,n})=
\big(\nabla_x+\frac{1}{2}\sum_{m=1}^{3}\Mb_m
x\frac{\partial}{\partial z_{m}}\big),$$ where
$\nabla_{x}=\big(\frac{\partial}{\partial
x_{11}},\ldots,\frac{\partial}{\partial x_{4\,n}}\big)$ the {\it
horizontal gradient}. Any vector field $Y$ belonging to $\mathcal
T Q^n$ is called the {\it horizontal vector field}. In particular,
the horizontal gradient $X$ is a horizontal vector field, that
justifies the name ``horizontal" gradient.

The commutation relations are as follows
$$
[X_{1l},X_{2l}]=-a_{1l}Z_{1},\quad
[X_{1l},X_{3l}]=a_{3l}Z_{3},\quad [X_{1l},X_{4l}]=a_{2l}Z_{2},
$$ $$
[X_{2l},X_{3l}]=a_{2l}Z_{2}, \quad
[X_{2l},X_{4l}]=-a_{3l}Z_{3},\quad [X_{3l},X_{4l}]=-a_{1l}Z_{1},$$
for any $l=1,\ldots,n$.

A basis of one-forms dual to $X_{kl},Z_{m}$, is given by
$dx_{kl},\vartheta_{m}$, with
\begin{equation}\label{74}
\vartheta_{m} = dz_{m}-\frac{1}{2}\big(\mathbf M_m x,dx\big).
\end{equation}
Since the interior product $\vartheta_{m}(X_{kl})$ vanishes for
all $m=1,2,3$, $k=1,\ldots,4$, $l=1,\ldots,n$ we have the product
$\vartheta_{m}(Y)$ vanishing on all horizontal vector fields $Y$.

\begin{remark}\label{r1} If we formally put $a_{1l}=0$, $l=1,\ldots,n$, then we obtain another
example of a quaternion anisotropic group with $2$-dimensional
center. The case $a_{1l}=a_{2l}=0$, $l=1,\ldots,n$, corresponds to
an anisotropic group with $1$-dimensional center.
\end{remark}

\section{Horizontal curves and their geometric characteristics}

Summarizing the results of the previous section we can say that
$Q^n$ is a space of ($4n+3$)~-tuples of real numbers $\mathbb
R^{4n+3}$ where the commutative group operation ``$+$" is replaced
by the noncommutative law ``$\circ$". Respectively, the left
translation $L_x(x^{\prime})=x+x^{\prime}$ (that in the commutative
case coincides with the right translation) is substituted by the
left translation $L_q(q^{\prime})=q\circ q^{\prime}$. The
corresponding Lie algebras are fundamentally different. The constant
vector fields $X_i=\frac{\partial}{\partial x_i}$, $i=1,\ldots,
4n+3$, of Euclidean space are replaced by the vector
fields~\eqref{62} and~\eqref{62z}. Moreover, since the group of
vector fields~\eqref{62z} is completely generated by the group
~\eqref{62} by means of commutation relations, the geometry of the
group $Q^n$ is defined by the horizontal bundle $\T Q^n$. The
velocity and the distance should respect the horizontal bundle $\T
Q^n$. Since $[X_{k_il_j},X_{k_al_b}]\notin \T Q^n$ the horizontal
bundle is not integrable, i.e., there is no surface locally tangent
to it~\cite{AMacK}. As we say, the geometry is defined by the
horizontal bundle, so it is sufficient to define the Riemanian
metric only on the horizontal bundle $\T Q^n$ of tangent bundle
of~$Q^n$. These kinds of manifolds have aquired the name
subRiemanian manifolds. The definitions and basic notations of
subRiemanian geometry can be found in,~e.g.,~\cite{Str}.

A continuous map $c(s):[0,1]\to Q^n$ is called a curve. We say
that a curve $c(s)$ is {\it horizontal} if its tangent vector
$\dot c(s)$ (if it exists) belongs to $\T Q^n$ at each point
$c(s)$. In other words, there are (measurable) functions
$a_{kl}(s)$ such that $\dot c(s)=\sum_{k,l}a_{kl}(s)X_{kl}(c(s))$.
We present some simple propositions that describe the geometry of
the group $Q^n$.

\begin{proposition}\label{p1} A curve
$c(s)=(x(s),z(s))$ is horizontal if and only if
\begin{equation}\label{1}
\dot z_{m} = \frac{1}{2}(\Mb_m x,\dot x),\quad m=1,2,3,
\end{equation} where $\dot x=(\dot x_{11},\dot x_{21},\dot x_{31},\dot x_{41},
\ldots,\dot x_{1n},\dot x_{2n},\dot x_{3n},\dot x_{4n})$.
\end{proposition}

\begin{proof}
We can write the tangent vector $\dot c(s)$ in the form
\begin{equation*}
\begin{split}
\dot c(s) = & (\dot x(s),\dot z(s)) = \sum_{k,l} \dot
x_{kl}(s)\frac{\partial}{\partial x_{kl}}+\sum_m\dot
z_{m}(s)\frac{\partial}{\partial z_{m}}
\\ = &
\Big(\dot x(s),\nabla_{x}+\frac{1}{2}\sum_{m}\Mb_m
x\frac{\partial}{\partial z_{m}}\Big)+\sum_m\Big(\dot
z_{m}(s)-\frac{1}{2}\big(\Mb_m x,\dot
x(s)\big)\Big)\frac{\partial}{\partial z_{m}}
\\ = &
\Big(\dot x(s),X(s)\Big)+\sum_m\Big(\dot
z_{m}(s)-\frac{1}{2}\big(\Mb_m x,\dot
x(s)\big)\Big)\frac{\partial}{\partial z_{m}}.
\end{split}
\end{equation*}
It is clear that $\dot c(s)$ is horizontal if and only if the
coefficients in front of $\frac{\partial}{\partial z_{m}}$,
$m=1,2,3$, vanish. This proves Proposition~\ref{p1}.
\end{proof}

\begin{corollary}
If a curve $c(s)=(x(s),z(s))$ is horizontal, then
$$\dot c(s)=\Big(\dot x(s),X\Big)=\sum_{k,l}\dot x_{kl}(s)X_{kl}.$$
\end{corollary}

It is easy to see the following statement.

\begin{proposition}\label{40}
The left translation $L_{q}$ of a horizontal curve
$c(s)=(x(s),z(s))$ is a horizontal curve $\tilde c(s)=L_{q}(c(s))$
with the velocity
\begin{equation}\label{16}\dot{\tilde c}(s)=(L_{q})_{*}\dot c(s)
=\sum_{k,l}\dot x_{kl}(s)X_{kl}(\tilde c(s))=\big(\dot
x(s),X(\tilde c(s))\big).\end{equation}
\end{proposition}

\begin{proposition}\label{39}
The acceleration vector $\ddot c(s)$ of a horizontal curve $c(s)$
is horizontal.
\end{proposition}

\begin{proof}
Let $c(s)$ be a horizontal curve. Then $\dot c(s)\in \T Q^n_{c(s)}$.
Let us show that $\ddot c(s)\in \T Q^n_{c(s)}$. Differentiating
equalities~\eqref{1} of the horizontality condition and making use
of $(\Mb_mx,x)=0$ for $m=1,2,3$ and any $x\in \mathbb R^{4n}$, we
deduce that
$$\ddot z_{m}(s)=\frac{1}{2}\Big(\big(\Mb_m\dot x(s),\dot x(s)\big)
+\big(\Mb_m x(s),\ddot x(s)\big)\Big) =\frac{1}{2}\big(\Mb_m
x(s),\ddot x(s)\big)$$ for $m=1,2,3$. Then the acceleration vector
along $c(s)$ is
\begin{equation*}
\begin{split}
\ddot c(s) = &\big(\ddot x(s),\nabla_x\big)+\big(\ddot
z(s),\nabla_z\big) = \Big(\ddot x(s),\nabla_{x}+ \frac{1}{2}\sum_m
\Mb_m x(s) \frac{\partial}{\partial
z_{m}}\Big) \\
+ & \sum_m\Big(\ddot z_{m}(s)-\frac{1}{2}\big(\ddot x(s),\Mb_m
x(s)\big)\Big)\frac{\partial}{\partial z_{m}}  = \big(\ddot
x(s),X(c(s))\big).
\end{split}
\end{equation*}
This means that the vector $\ddot c(s)$ is horizontal. The
proposition is proved.
\end{proof}

\section{Hamiltonian formalism}

In this section we study the geometry of the anisotropic quaternion
group $Q^n$ making use of the Hamiltonian formalism. The geometry of
the group is induced by the sub-Laplacian
$\Delta_0=\sum_{k,l}X^{2}_{kl}$. Operators of such type are studied,
for instance, in~\cite{BGGR3,CDG}. Since the vector fields $X_{kl}$
satisfying the Chow's condition, by a theorem of
H\"ormander~\cite{Hor}, the operator $\Delta_0$ is hypoelliptic.
Explicitly, the sub-Laplacian has the form:
\begin{equation*}
\Delta_0=\sum_{l=1}^{n}\sum_{k=1}^{4}X^{2}_{kl} = \Big(
\Delta_{x}+\frac{1}{4}\sum_{m=1}^{3}\Big(\sum_{l=1}^{n}a^2_{ml}|x_l|^2\Big)\frac{\partial^2}{\partial
z_m^2}+\sum_{m=1}^{3}\big(\Mb_m x,
\nabla_{x}\big)\frac{\partial}{\partial z_{m}}\Big),
\end{equation*} where $\nabla_{x}=(\frac{\partial}{\partial
x_{11}},\ldots,\frac{\partial}{\partial x_{4n}})$,
$\Delta_{x}=\sum_{l=1}^{n}\sum_{k=1}^{4}\frac{\partial^2}{\partial
x_{kl}^2}$, and $|x_l|^2=\sum_{k=1}^{4}x_{kl}^2$. To present the
Hamiltonian function we introduce the formal variables
$\xi=(\xi_{11},\ldots,\xi_{4n})$ with
$\xi_{kl}=\frac{\partial}{\partial x_{kl}}$ and
$\theta=(\theta_1,\theta_2,\theta_3)$ with
$\theta_m=\frac{\partial}{\partial z_{m}}$, $m=1,2,3$. The
associated with sub-Laplacian $\Delta_0$ Hamiltonian function
$H(\xi,\theta,x,z)$ is following
\begin{equation}
H(\xi,\theta,x,z) =
|\xi|^2+\frac{1}{4}\sum_{m=1}^{3}\Big(\sum_{l=1}^{n}a_{ml}^2|x_l|^2\Big)\theta_m^2
+ \Big(\big(\sum_{m=1}^{3}\theta_{m}\Mb_m\big) x,\xi\Big),
\end{equation}
where $|\xi|^2=\sum_{l,k}\xi_{kl}^2$, and diagonal blocks of the
matrix $\sum_{m=1}^{3}\theta_{m}\Mb_m$ are of the form
$$\left[\array{rrrr} 0 & a_{1l}\theta_{1} & -a_{3l}\theta_{3} & -a_{2l}\theta_{2}
\\
-a_{1l}\theta_{1} & 0 & -a_{2l}\theta_{2} & a_{3l}\theta_{3}
\\
a_{3l}\theta_{3} & a_{2l}\theta_{2} & 0 & a_{1l}\theta_{1}
\\
a_{2l}\theta_{2} & -a_{3l}\theta_{3} & -a_{1l}\theta_{1} &
0\endarray\right].$$ We use the following notation:
$$
A_m^2=\left[\array{ccc}a_{m1}^2\U & \ & 0
\\ \ & \ddots & \
\\
 0  & \ &a_{mn}^2\U
\endarray\right],
$$
$\Tb^2=\sum_{m=1}^{3}\theta_{m}^2A_m^2$, and
$\Mb=\sum_{m=1}^{3}\theta_{m}\Mb_m$. We also introduce some
different metrics for convenience:
$|\theta|^2_l=\sum_{m=1}^{3}\theta^2_ma^2_{ml}$,
$|x|^2_{B}=(B^2x,x)=(Bx,Bx)$, where $B$ is a diagonal matrix. In
this notation we get $\sum_{l=1}^{n}a_{ml}^2|x_l|^2=|x|^2_{A_m}$.
The Hamiltonian function takes a new form in this notation
\begin{equation}\label{35}
H(\xi,\theta,x,z) =
|\xi|^2+\frac{1}{4}\sum_{m=1}^{3}|x|^2_{A_m}\theta_m^2 + \big(\Mb
x,\xi\big)=|\xi|^2+\frac{1}{4}(\Tb^2 x,x) + \big(\Mb x,\xi\big),
\end{equation}
and the corresponding Hamiltonian system obtains the form
\begin{equation}\label{2}
\begin{split} \left\{\array{lll}\dot x &= \frac{\partial H}{\partial
\xi}=2\xi+\Mb x
\\
\dot z_{m} &=  \frac{\partial H}{\partial
\theta_{m}}=\frac{\theta_{m}}{2}|x|^2_{A_m}+\big(\Mb_m
x,\xi\big),\quad m=1,2,3.
\\
\dot \xi &= -\frac{\partial H}{\partial x}=-\frac{1}{2}\Tb^2
x+\Mb\xi
\\
\dot \theta_m &= -\frac{\partial H}{\partial z_m}=0.
\endarray\right.
\end{split} \end{equation}
The solutions $\gamma(s)=(x(s),z(s),\xi(s),\theta(s))$ of the
system~\eqref{2} are called {\it bicharacteristics}.

\begin{definition}
Let $P_1(x_0,z_0)$, $P_2(x,z)\in Q^n$. A geodesic from $P_1$ to
$P_2$ is the projection of a bicharacteristic $\gamma(s)$,
$s\in[0,\tau]$, onto the $(x,z)$-space, that satisfies the
boundary conditions
$$\big(x(0),z(0)\big)=(x_0,z_0),\qquad \big(x(\tau),z(\tau)\big)=(x,z).$$
\end{definition}

The next properties of the matrices $\M_m$, $m=1,2,3$, are
obvious: \begin{equation}\label{m1} \M_m^2=-\mathcal U,\ \
m=1,2,3,\ \text{where}\ \mathcal U\ \text{is the unit}\ (4\times
4)-\text{matrix.}\end{equation}
\begin{eqnarray}\label{m2} \M_1\M_2=-\M_2\M_1 & = & \M_3,\quad
\M_2\M_3=-\M_3\M_2=\M_1,\nonumber\\
\M_3\M_1 & = & -\M_1\M_3=\M_2.\end{eqnarray}
\begin{equation}\label{m3} \M_m^{-1}=-\M_m,\ \text{where}\  \M_m^{-1}\ \text{is the inverse
matrix of}\  \M_m,\ m=1,2,3.\end{equation}
\begin{equation}\label{m4}\M_m^{T}=-\M_m, \ \text{where}\  \M^{T}_m\ \text{is the
transposed matrix for}\  \M_m,\ m=1,2,3.\end{equation}
\begin{equation}\label{m5}(\M_m x,x)=0,\ m=1,2,3,\ \text{for
any}\  x\in\mathbb R^{4}.\end{equation} As a corollary we obtain
some useful formulas.

\begin{proposition} In the above-mentioned notations we have
\begin{equation}\label{an1}(\Mb_mx,\Mb x)=\theta_{m}\sum_{l=1}^{n}a_{ml}^2|x_l|^2=\theta_m|x|^2_{A_m}\quad \text{for
any}\quad j=1,2,3.\end{equation}
\begin{equation}\label{an2}\Mb^2=-\sum_{m=1}^{3}\theta_{m}^2A_m^2=-\Tb^2,\quad \Mb^3=-\Tb^2\Mb,\quad
\Mb^4=\Tb^4,\quad \Mb^5=\Tb^4\Mb\ \ldots\end{equation}
\end{proposition}
\begin{proof}
\begin{eqnarray*}(\Mb_mx,\Mb
x) & = &
(\Mb_mx,\sum_{j=1}^{3}\theta_j\Mb_jx)=\sum_{j=1}^{3}\theta_j(\Mb_mx,\Mb_j
x) \\ & = &
\sum_{j=1}^{3}\theta_j\Big((a_{j1}\M_jx_1,\ldots,a_{jn}\M_jx_n)\
,\ (a_{m1}\M_mx_1,\ldots,a_{mn}\M_mx_n)\Big)
\\ & = & \sum_{j=1}^{3}\theta_j\sum_{l=1}^{n}a_{jl}a_{ml}(\M_jx_l,\M_mx_l)
=\sum_{j=1}^{3}\theta_j\sum_{l=1}^{n}a_{jl}a_{ml}(-\M_m\M_jx_l,x_l)
\\ & = & \theta_m\sum_{l=1}^{n}a_{ml}^2|x_l|^2=\theta_m|x|^2_{A_m}
\end{eqnarray*} by the properties~\eqref{m1},~\eqref{m2},~\eqref{m4}, and~\eqref{m5} of matrix
$\M_m$.

To prove~\eqref{an2} we note that $\Mb_m\Mb_j=-\Mb_j\Mb_m$
by~\eqref{m2} and $\Mb_m^2=-A^2_m$ by the property~\eqref{m1} for
any $j,m=1,2,3$. Then
$$\Mb^2=\sum_{j=1}^{3}\sum_{m=1}^{3}\theta_j\theta_m\Mb_j\Mb_m=\sum_{m=1}^{3}\theta_m^2\Mb_m^2=
-\sum_{m=1}^{3}\theta_m^2A^2_m=-\Tb^2.$$ The rest is obvious.
\end{proof}

\begin{lemma}\label{p2}
Any geodesic is a horizontal curve.
\end{lemma}

\begin{proof}
Let $c(s)=(x(s),z(s))$ be a geodesic. The system~\eqref{2} implies
\begin{equation}\label{17}
\dot z_{m}=\frac{\theta_{m}}{2}|x|^2_{A_m}+\frac{1}{2}\big(\Mb_m
x,2\xi\big)=\frac{\theta_{m}}{2}|x|^2_{A_m}+\frac{1}{2}\big(\Mb_m
x,\dot x\big)+\frac{1}{2}\big(\Mb_m x,2\xi-\dot x\big).
\end{equation} Making use of the first line
of the system~\eqref{2}, we write the last term of~\eqref{17} as
\begin{equation}\label{18}
\frac{1}{2}\big(\Mb_m x,2\xi-\dot x\big) = -\frac{1}{2}(\Mb_m
x,\Mb x)=-\frac{\theta_{m}}{2}|x|^2_{A_m}.
\end{equation} Here we used the formula~\eqref{an1}.
Combining~\eqref{17} and~\eqref{18} we deduce
\begin{equation}\label{an3}\dot
z_{m}=\frac{\theta_{m}}{2}|x|^2_{A_m}+\big(\Mb_m
x,\xi\big)=\frac{1}{2}\big(\Mb_m x,\dot x\big),\quad
m=1,2,3.\end{equation} Therefore, $c(s)$ is a horizontal curve by
Proposition~\ref{p1}.
\end{proof} Lemma~\ref{p2} shows that the second equation of the
system~\eqref{2} is nothing more then the horizontality
condition~\eqref{1}.

Let us try to solve the Hamiltonian system explicitly. The last
equation in~\eqref{2} shows that the function $H(\xi,\theta,x,z)$
does not depend on $z$. We obtain that $\theta_{m}$ are constants
which can be used as Lagrangian multipliers. Multiplying the first
line of system~\eqref{2} by $\Mb$, we obtain
\begin{equation}\label{5}\Mb\dot x=2\Mb\xi-\Tb^2 x.\end{equation}
Expressing $\Mb\xi$ from~\eqref{5} and substituting it in the
equation for $\dot \xi$ from~\eqref{2}, we get
\begin{equation}\label{20}\dot\xi=\frac{\Mb\dot
x}{2}.\end{equation} We differentiate the first equation of
~\eqref{2} and substitute the $\dot\xi$ from~\eqref{20}. Finally,
we deduce \begin{equation}\label{9}\ddot x=2\dot\xi+\Mb\dot
x=2\Mb\dot x.\end{equation}

Let us solve the equation~\eqref{9}. We substitute $y(s)=\dot
x(s)$. The equation $\dot y(s)=2\Mb y(s)$ has a solution
$y(s)=\exp(2s\Mb)y(0)$. Therefore,
\begin{equation}\label{h1}\dot x(s)=\exp(2s\Mb)\dot x(0).\end{equation} Let us discuss the properties of the
matrix $\exp(2s\Mb)$. For simplicity of notation we write $[B]_l$
for $l$-block of a block diagonal matrix $B$.

\begin{lemma}\label{49}
The exponent $\exp\big(2s\Mb\big)$ is an antisymmetric block matrix
that commutes with $\Mb$ and which blocks can be written in the
form:
\begin{equation}\label{exp}[\exp(2s\Mb)]_l=\cos(2s|\theta|_l)\U
+\frac{\sin(2s|\theta|_l)}{|\theta|_l}[\Mb]_l.\end{equation}
\end{lemma}
\begin{proof}
We observe that
\begin{equation*}
\begin{split}
[\exp\big(2s\Mb)]_l = &
\sum^{\infty}_{n=0}\frac{(2s)^n}{n!}[\Mb]_l^n
 = \mathbf U\sum^{\infty}_{k=0}\frac{(2s|\theta|_l)^{4k}}{(4k)!}+
\frac{[\Mb]_l}{|\theta|_l}\sum^{\infty}_{k=0}\frac{(2s|\theta|_l)^{4k+1}}{(4k+1)!}
\\ &-
\mathbf U\sum^{\infty}_{k=0}\frac{(2s|\theta|_l)^{4k+2}}{(4k+2)!}-
\frac{[\Mb]_l}{|\theta|_l}\sum^{\infty}_{k=0}\frac{(2s|\theta|_l)^{4k+3}}{(4k+3)!}
\end{split}
\end{equation*} by~\eqref{an2}. We conclude that the matrices
$\Mb$ and $\exp \big(2s\Mb\big)$ commute. Note that
$$\sum^{\infty}_{k=0}\frac{(2s|\theta|_l)^{4k}}{(4k)!}-\sum^{\infty}_{k=0}\frac{(2s|\theta|_l)^{4k+2}}{(4k+2)!}
=\cos(2s|\theta|_l)$$ and
$$\sum^{\infty}_{k=0}\frac{(2s|\theta|_l)^{4k+1}}{(4k+1)!}-
\sum^{\infty}_{k=0}\frac{(2s|\theta|_l)^{4k+3}}{(4k+3)!}=\sin(2s|\theta|_l).$$
With this notation, $\exp\big(2s\Mb\big)$ is a block diagonal
matrix with the blocks $[\exp\big(2s\Mb\big)]_l =$
\begin{eqnarray*}
= \left[\array{llll} \cos(2s|\theta|_l) &
\frac{a_{1l}\theta_{1}}{|\theta|_l}\sin(2s|\theta|_l) &
-\frac{a_{3l}\theta_{3}}{|\theta|_l}\sin(2s|\theta|_l) &
-\frac{a_{2l}\theta_{2}}{|\theta|_l}\sin(2s|\theta|_l)
\\
-\frac{a_{1l}\theta_{1}}{|\theta|_l}\sin(2s|\theta|_l) &
\cos(2s|\theta|_l) &
-\frac{a_{2l}\theta_{2}}{|\theta|_l}\sin(2s|\theta|_l) &
\frac{a_{3l}\theta_{3}}{|\theta|_l}\sin(2s|\theta|_l)
\\
\frac{a_{3l}\theta_{3}}{|\theta|_l}\sin(2s|\theta|_l) &
\frac{a_{2l}\theta_{2}}{|\theta|_l}\sin(2s|\theta|_l) &
\cos(2s|\theta|_l) &
\frac{a_{1l}\theta_{1}}{|\theta|_l}\sin(2s|\theta|_l)
\\
\frac{a_{2l}\theta_{2}}{|\theta|_l}\sin(2s|\theta|_l) &
-\frac{a_{3l}\theta_{3}}{|\theta|_l}\sin(2s|\theta|_l) &
-\frac{a_{1l}\theta_{1}}{|\theta|_l}\sin(2s|\theta|_l) &
\cos(2s|\theta|_l)
\endarray\right].
\end{eqnarray*}
\end{proof}

The group structure allows to restrict our considerations to the
curves issuing from the origin. Hence, $x(0)=0$. The
equation~\eqref{h1} has the form
\begin{equation}\label{47} \dot x_l(s)=\cos(2s|\theta|_l)\U\dot
x_l(0)+\frac{\sin(2s|\theta|_l)}{|\theta|_l}[\Mb]_l\dot
x_l(0),\quad l=1,\ldots,n,
\end{equation} by~\eqref{exp}.
Integrating from $0$ to $s$ we get
\begin{equation}\label{48} x_l(s)=\frac{1-\cos(2s|\theta|_l)}{2|\theta|_l^2}[\Mb]_l\dot x_l(0)
+\frac{\sin(2s|\theta|_l)}{2|\theta|_l}\U\dot x_l(0),\quad
l=1,\ldots,n .
\end{equation}
Let us describe the $z$-components of a geodesic curve. If a curve
is geodesic, then it is horizontal by Lemma~\ref{p2}, and we have
\begin{equation*}
\begin{split}
\dot z_{m}(s) = \frac{1}{2}\big(\Mb_m x(s),\dot x(s)\big)= &
\sum_{l=1}^{n}\Bigl(\frac{\cos(2s|\theta|_l)(1-\cos(2s|\theta|_l))}{4|\theta|_l^2}\big([\Mb_m]_l[\Mb]_l\dot
x_l(0),\dot x_l(0)\big)
\\ & +  \frac{\sin(2s|\theta|_l)(1-\cos(2s|\theta|_l))}{4|\theta|_l^3}\big([\Mb_m]_l[\Mb]_l\dot
x_l(0),[\Mb]_l\dot x_l(0)\big)\\ & +
\frac{\sin(2s|\theta|_l)\cos(2s|\theta|_l)}{4|\theta|_l}\big([\Mb_m]_l\dot
x_l(0),\dot x_l(0)\big)\\ & +
\frac{\sin^2(2s|\theta|_l)}{4|\theta|_l^2}\big([\Mb_m]_l\dot
x_l(0),[\Mb]_l\dot x_l(0)\big)\Big), \end{split}\end{equation*}
for $m=1,2,3$ by~\eqref{47} and~\eqref{48}. The
properties~\eqref{m5} and~\eqref{an1} imply
$$\big([\Mb_m]_l\dot x_l(0),\dot x_l(0)\big)=([\Mb_m]_l[\Mb]_l\dot x_l(0),[\Mb]_l\dot
x_l(0)\big)=0,$$ and
$$\big([\Mb_m]_l[\Mb]_l\dot x_l(0),\dot x_l(0)\big)=-\big([\Mb_m]_l\dot x_l(0),[\Mb]_l\dot
x_l(0)\big)=-\theta_{m}a_{ml}^2|\dot x_l(0)|^2.$$ Finally, we see
that
\begin{equation}\label{50}\dot z_{m}(s)
= \sum_{l=1}^{n}\Bigl(\frac{\theta_{m}a^2_{ml}|\dot
x_l(0)|^2}{4|\theta|_l^2}\big(1-\cos(2s|\theta|_l)\big)\Bigr),\quad
m=1,2,3.
\end{equation}
Integrating equations~\eqref{50}, we get
\begin{equation}\label{53}
z_{m}(s) = \sum_{l=1}^{n}\Bigl(\frac{\theta_{m}a^2_{ml}|\dot
x_l(0)|^2}{4|\theta|_l^2}\big(s-\frac{\sin(2s|\theta|_l)}{2|\theta|_l}\big)\Bigr),
\quad m=1,2,3.
\end{equation}

\begin{lemma} Not all of horizontal curves are geodesics.
\end{lemma}

\begin{proof}
To prove this proposition we present an example. The curve
$$c(s)=\Big(\frac{s^2}{2},s,\frac{s^2}{2},s,0,\ldots,0,\frac{a_{11}s^3}{6},c_1,c_2\Big)$$  is
horizontal with $c_1,c_2$ constant. Indeed, \begin{equation*}
\begin{split}\dot
z_{1}(s) = & \frac{a_{11}s^2}{2}, \quad \frac{1}{2}\big(\Mb_1
x,\dot x\big) =
\frac{a_{11}}{2}(s^2-\frac{s^2}{2}+s^2-\frac{s^2}{2})=\frac{a_{11}s^2}{2},\\
\dot z_{2}(s)= & 0, \quad \frac{1}{2}\big(\Mb_2 x,\dot
x\big)=  \frac{a_{21}}{2}(-s^2-\frac{s^2}{2}+s^2+\frac{s^2}{2})=0,\\
\dot z_{3}(s)= & 0, \quad \frac{1}{2}\big(\Mb_3 x,\dot x\big)=
\frac{a_{31}}{2}(-\frac{s^3}{2}+s+\frac{s^3}{2}-s)=0.\end{split}
\end{equation*} From
the other hand, the curve $c(s)$ does not satisfy the
system~\eqref{9}. The system~\eqref{9} gets the form
\begin{equation*}
\begin{split}
\left\{\array{lll}1 =
2(a_{11}\theta_{1}-a_{31}\theta_{3}s-a_{21}\theta_{2})
\\
0 = 2(-a_{11}\theta_{1}s-a_{21}\theta_{2}s+a_{31}\theta_{3})
\\
1 = 2(a_{31}\theta_{3}s+a_{21}\theta_{2}+a_{11}\theta_{1}) \\ 0 =
2(a_{21}\theta_{2}s-a_{31}\theta_{3}-a_{11}\theta_{1}s)\endarray
\right.
\end{split}
\end{equation*} for the curve $c(s)$. Summing up the first and
the third equation, and then, the second and the forth ones, we
write the latter system as follows
\begin{equation*}
\begin{split}
\left\{\array{lll}2 = 4a_{11}\theta_{1}
\\
0 = -4a_{11}\theta_{1}s
\\
1 = 2(a_{31}\theta_{3}s+a_{21}\theta_{2}+a_{11}\theta_{1})
\\
0 =
2(a_{21}\theta_{2}s-a_{31}\theta_{3}-a_{11}\theta_{1}s).\endarray
\right.
\end{split}
\end{equation*} We see that the first and the second equations
contradict each other.
\end{proof}

\begin{lemma}\label{28}
A curve $c$ is a geodesic for the group $Q^n$ if and only if
\begin{itemize}\item[(i)]{$c(s)$ is a horizontal curve
and}\item[(ii)]{$c(s)$ satisfies $\ddot c(s)=2\Mb\dot c(s)$,
$l=1,\ldots,n$.}\end{itemize}
\end{lemma}
\begin{proof}
If a curve is geodesic, then  it is horizontal by Lemma~\ref{p2}.
Proposition~\ref{39} implies that the vector $\ddot c$ is also
horizontal: $\ddot c=\sum_{l=1}^{n}\ddot x_lX_l$. Since $\ddot
x(s)=2\Mb\dot x$ by~\eqref{9}, we obtain the necessary result.

Let the curve $c(s)$ satisfy~(i) and (ii) of Lemma~\ref{28}. The
horizontality condition~(i) of Lemma~\ref{28} can be written in
the form \begin{equation}\label{21}\dot
z_{m}=\frac{1}{2}\Big(\Mb_m x,\dot
x\Big)=\frac{\theta_{m}}{2}|x|^2_{A_m}+\big(\Mb_m
x,\xi\big)=\frac{\partial H}{\partial\theta_{m}},\quad
m=1,2,3.\end{equation} as in~\eqref{an3}. We see that $c(s)$
satisfies the equations of the second line of~\eqref{2}. The
condition~(ii) of Lemma~\ref{28} admits the form $\ddot
x(s)=2\Mb\dot x(s)$ in the coordinate functions. Define the
following curve $\gamma(s)=(x(s),z(s),\xi(s),\theta)$ in the
cotangent space, where
\begin{equation}\label{22}\xi=\frac{\dot x(s)}{2}-\frac{1}{2}\Mb
x(s)\quad\text{with}\quad
\theta=(\theta_1,\theta_2,\theta_3)\quad\text{constant}.\end{equation}
The relations~\eqref{22} imply the equations of the first and the
last lines of~\eqref{2}. Differentiating~\eqref{22}, we get
$$\dot\xi=\frac{\ddot x}{2}-\frac{1}{2}\Mb\dot x=\Mb\dot
x-\frac{\Mb\dot x}{2}=\frac{1}{2}\Mb(2\xi+\Mb
x)=\Mb\xi-\frac{1}{2}\Tb^2 x,$$ by the condition~(ii) of
Lemma~\ref{28}, \eqref{22}, and~\eqref{an2}. Thus, $\gamma(s)$
satisfies the Hamilton system~\eqref{2}. Then, the projection onto
the $(x,z)$-space, that coincides with $c(s)$, is a geodesic.
\end{proof}

\section{Connectivity by geodesics.}

Let us ask in the following question. Is it possible to join
arbitrary two points of $Q^n$ by a horizontal curve? A theorem by
Chow~\cite{Ch} gives an affirmative answer. We present a direct
proof and calculate the number of geodesics connecting the origin
with different points. We need the next simple observation.

\begin{proposition}\label{55} The kinetic energies $\mathcal E=\frac{1}{2}|\dot x|^2$,
$\mathcal E_m=\frac{1}{2}|\dot x|^2_{A_m}$ are preserved along
geodesics.
\end{proposition}

\begin{proof} In fact,
\begin{equation*} \frac{d\mathcal E_m}{ds}=\big(A_m\dot
x,A_m\ddot x\big)=2\big(A_m\dot x, \Mb A_m\dot x\big)=0
\end{equation*} by Lemma~\ref{28} and property~\eqref{m5} of the matrices $\M_m$.
The same is for $\mathcal E$.
\end{proof}

\subsection{Connectivity between $(0,0)$ and $(x,0)$, $x\neq 0$.}

\begin{theorem}\label{54}
A smooth curve $c(s)$ is horizontal with constant $z$-coordinates
$z_1,z_2,z_3$ if and only if
$c(s)=(\alpha_{11}s,\ldots,\alpha_{4\,n}s,z_1,z_2,z_3)$ with
$\alpha_{kl}\in\mathbb R$ and
$\sum_{l=1}^{n}\sum_{k=1}^{4}\alpha_{kl}^2\neq 0$. In other words,
there is only one geodesic joining the origin with a point
$(x,0)$.
\end{theorem}
\begin{proof}
Let $c(s)$ be a horizontal curve with constant $z$-coordinates
$z_1,z_2,z_3$. Then $\dot z_m=0$ and~\eqref{50} implies
$$
0=\dot z_{m} = \theta_{m}\sum_{l=1}^{n}\frac{a^2_{ml}|\dot
x_l(0)|^2}{2|\theta|_l^2}\sin^2\big( s|\theta|_l\big),\quad
m=1,2,3.
$$ We define by continuity $\frac{\sin^2
\big(s|\theta|_l\big)}{|\theta|_l^2}=s^2$ at $|\theta|_l=0$. Since
the sum $\sum_{l=1}^{n}\frac{a^2_{ml}|\dot
x_l(0)|^2}{2|\theta|_l^2}\sin^2 \big(s|\theta|_l\big)$ is not
identically zero we deduce, that $\theta_m=0$ for all $m=1,2,3$.
The Hamiltonian system~\eqref{2} is reduced to the next one
\begin{equation*}
\begin{split} \left\{\array{lll}\dot x &= 2\xi
\\
0 &= \big(\Mb_m x,\xi\big),\quad m=1,2,3
\\
\dot \xi &= 0
\\
\theta_m &= 0.
\endarray\right.
\end{split} \end{equation*}
We see that $\xi$ is a constant vector. Taking into account that
$x(0)=0$, we get $x(s)=(\alpha_{11}s,\ldots,\alpha_{4\,n}s)$ with
$\alpha_{kl}=2\xi_{kl}$. This proves the statement.

Now, let us assume that
$c(s)=(\alpha_{11}s,\ldots,\alpha_{4\,n}s,z_1,z_2,z_3)$ with
constant $z$-components. Set $\alpha
s=(\alpha_{11}s,\ldots,\alpha_{4\,n}s)$. Recall, that $(\Mb_m
\alpha,\alpha)=0$ for any vector
$\alpha=(\alpha_{11},\ldots,\alpha_{4\,n})$ and $m=1,2,3$. Then,
$$\dot z_{m}=  0  = \frac{1}{2}(\Mb_m(\alpha s),\dot{(\alpha s)})=\frac{s}{2}(\Mb_m \alpha,\alpha),
\qquad  m=1,2,3.$$ The horizontal condition~\eqref{1} holds for
all three $z$-components.
\end{proof}

\subsection{Connectivity between $(0,0)$ and $(0,z)$, $z\neq 0$.}

We need to solve the equation~\eqref{9} with the boundary
conditions
$$x(0)=x(1)=z(0)=0,\qquad z(1)=z.$$ We also need to know the
initial velocity $\dot x(0)$ since we do not have enough
information about the behavior of $x$-coordinates. In the
following theorem we use the notations $\nb=(n_1,\ldots,n_n)$,
$n_l\in\mathbb N$, for $l=1,\ldots, n$, and
\begin{equation}\label{ann} N= \left[\array{rrr} \frac{1}{\pi n_1}\mathcal U & \
 & 0
\\
\ & \ddots & \
\\
0 & \ & \frac{1}{\pi n_n}\mathcal U
\endarray\right].
\end{equation}

\begin{theorem}\label{4}
There are infinitely many geodesics joining the origin with a
point $(0,z)$. The corresponding equations for each
$\nb=(n_1,\ldots,n_n)$ are
\begin{eqnarray}\label{65} x^{(\nb)}_l(s) = 2\frac{1-\cos(2s\pi n_l)}{(\pi
n_l)^2}[\mathbf Z]_l\dot x_l(0) +\frac{\sin(2s\pi n_l)}{2\pi
n_l}\U\dot x_l(0),\ \ l=1,\ldots,n,\end{eqnarray} where $\mathbf
Z$ is a block diagonal matrix with the blocks
\begin{equation}\label{71} [\mathbf Z]_l= \left[\array{cccc} 0 &
\frac{z_{1}a_{1l}}{|\dot x(0)|^2_{NA_1}} &
-\frac{z_{3}a_{3l}}{|\dot x(0)|^2_{NA_3}} &
-\frac{z_{2}a_{2l}}{|\dot x(0)|^2_{NA_2}}
\\
-\frac{z_{1}a_{1l}}{|\dot x(0)|^2_{NA_1}} & 0 &
-\frac{z_{2}a_{2l}}{|\dot x(0)|^2_{NA_2}} &
\frac{z_{3}a_{3l}}{|\dot x(0)|^2_{NA_3}}
\\
\frac{z_{3}a_{3l}}{|\dot x(0)|^2_{NA_3}} &
\frac{z_{2}a_{2l}}{|\dot x(0)|^2_{NA_2}} & 0 &
\frac{z_{1}a_{1l}}{|\dot x(0)|^2_{NA_1}}
\\
\frac{z_{2}a_{2l}}{|\dot x(0)|^2_{NA_2}} &
-\frac{z_{3}a_{3l}}{|\dot x(0)|^2_{NA_3}} &
-\frac{z_{1}a_{1l}}{|\dot x(0)|^2_{NA_1}} & 0\endarray\right],
\end{equation} and
\begin{equation}\label{68} z^{(\nb)}_m(s) =
\frac{z_m}{|\dot x(0)|^2_{NA_m}}\sum_{l=1}^{n}\frac{a^2_{ml}|\dot
x(0)|^2}{|\pi n_l|^2}\Big(s-\frac{\sin(2s\pi n_l)}{2\pi n_l}\Big),
\qquad m=1,2,3.\end{equation} The lengths of corresponding
geodesics are
$$l_{\nb}^2  =  16\sum_{m=1}^{3}\frac{z_m^2(1)}{\sum_{l=1}^{n}\frac{a^2_{ml}|\dot
x_l(0)|^2}{(\pi n_l)^2}}=16\sum_{m=1}^{3}\frac{z_m^2(1)}{|\dot
x_l(0)|_{NA_m}^2}.$$
\end{theorem}
\begin{proof}
Substituting $s=1$ in~\eqref{48}, we calculate
\begin{eqnarray*}0 & = & |x(1)|^2=\sum_{l=1}^{n}|x_l(1)|^2=\sum_{l=1}^{n}\Big(\frac{(\cos
2|\theta|_l-1)^2}{4|\theta|_l^4}\big([\Mb]_l\dot
x_l(0),[\Mb]_l\dot x_l(0)\big)
+\frac{\sin^22|\theta|_l}{4|\theta|_l^2}|\dot x_l(0)|^2\Big)
\\ & = & \sum_{l=1}^{n}\frac{\sin^2|\theta|_l}{|\theta|_l^2}|\dot x_l(0)|^2.\end{eqnarray*}
Since the kinetic energy $\mathcal E=\frac{|\dot x(0)|^2}{2}$ does
not vanish, there are indexes $l$, such that $|\dot x_l(0)|^2\neq
0$. We deduce in this case that
$$|\theta|_l=\sqrt{\theta_{1}^2a^2_{1l}+\theta_{2}^2a_{2l}^2+\theta_{3}^2a_{3l}^2}=\pi n_l,\qquad
n_l\in\mathbb N.$$ If $|\dot x_l(0)|$=0, then the corresponding
$|\theta|_l$ are arbitrary. Equalities~\eqref{53} give for $s=1$
\begin{equation}\label{64}
z_{m}(1)=\theta_{m}\sum_{l=1}^{n}\frac{a^2_{ml}|\dot
x_l(0)|^2}{4(\pi n_l)^2}=\frac{\theta_m}{4} |\dot
x(0)|^2_{NA_m},\qquad m=1,2,3.
\end{equation}
We find the unknown constants $\theta_m= \frac{4z_m(1)}{|\dot
x(0)|^2_{NA_m}}$. Substituting $\theta_m$
in~\eqref{48},~\eqref{53}, we obtain the equations~\eqref{65}
and~\eqref{68} for geodesics.

To calculate the length of geodesics, we observe that
$\sum_{m=1}^{3}\theta_m^2a^2_{ml}=\pi^2n_l^2$ and deduce
$$
\sum_{m=1}^{3}z_{m}(1)\theta_m=\frac{1}{4}\sum_{m=1}^{3}\sum_{l=1}^{n}\frac{\theta_{m}^2a^2_{ml}|\dot
x_l(0)|^2}{(\pi n_l)^2}=\frac{1}{4}\sum_{l=1}^{n}\frac{|\dot
x_l(0)|^2}{(\pi n_l)^2}\sum_{m=1}^{3}
\theta_m^2a_{ml}^2=\frac{1}{4}|\dot x(0)|^2=\frac{\mathcal E}{2}
$$
from~\eqref{64}. The lengths of geodesics are
\begin{eqnarray}\label{an7}l_{\nb}^2 & = & \Big(\int_{0}^{1}|\dot x(s)|\,ds\Big)^2=|\dot
x(0)|^2=4\sum_{m=1}^{3}z_{m}\theta_m=16\sum_{m=1}^{3}\frac{z_m^2(1)}{\sum_{l=1}^{n}\frac{a^2_{ml}|\dot
x_l(0)|^2}{(\pi n_l)^2}} \nonumber \\ & = &
16\sum_{m=1}^{3}\frac{z_m^2(1)}{|\dot x(0)|^2_{NA_m}},\qquad
\nb=(n_1,\ldots,n_n),\quad n_l\in\mathbb N.\end{eqnarray}
\end{proof}

\begin{remark} Let us discuss the cardinality of the set of geodesics.
In the general case in Theorem~\ref{4}, when $a_{ml}$ are
different we obtain the countably many geodesics connecting the
origin with a point $P=(0,z)$. If the multiindex
$\nb=(n_1,\ldots,n_n)$ increases, then the geodesic rotates more
frequently around straight line, connecting origin with $P=(0,z)$
approaching to the line and in the limit we obtain a limit curve
of Hausdorff dimension $2$. We present a graphic of three geodesic
curves~\eqref{65},\eqref{68} with the initial velocity $\dot
x(0)=(\dot x_{11},0,0,0,\dot x_{21},0,0,0,\ 0,\ldots,0)$ and the
end point $P=(0,z_1,0,0)$, where $\dot x_{11},\dot x_{21},z_1$ are
different from zero. The corresponding multiindexes are $\nb=1$,
$\nb=2$ and $\nb=5$.

%%%%%%%%%%%%%%%%%%%%%%%%FIGURE 1%%%%%%%%%%%%%%%%%%%%%%%%%%%%%%
\begin{figure}[ht]
\centering \scalebox{0.8}{\includegraphics{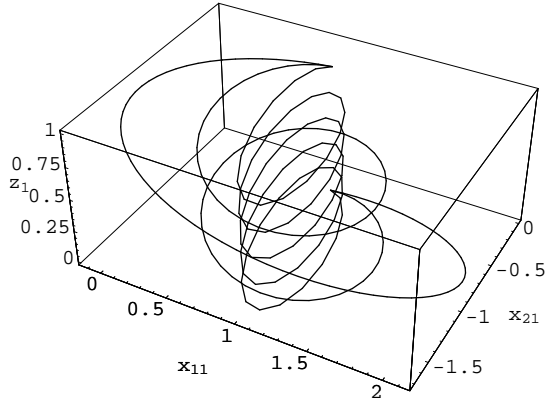}}
\caption[]{The graphs of geodesics where the vertical axis
represents $z$ coordinate and horizontal exes are $x_{11}$ and
$x_{21}$}\label{fig3}
\end{figure}
%%%%%%%%%%%%%%%%%%%%%%%%%%%%%%%%%%%%%%%%%%%%%%%%%%%%%%%%%%%%%%

The projection of geodesics into the horizontal subspace belongs
to an ellipsoid passing through the origin.

If, in particulary, $a_{1l}=a_{2l}=a_{3l}=a_l$ then rearranging the
indexes, we can assume that $a_1<a_2<\ldots<a_p=a_{p+1}=\ldots=a_n$.
Applying the rotation to the geodesics in the subspace
$(0,\ldots,0,x_{p},x_{p+1},\ldots,x_n, 0,0,0)$, we get uncountably
many geodesics. In this case we have the following estimate of their
lengths:
$$l_{\nb}^2=\frac{16|z(1)|^2}{\sum_{l=1}^{n}\frac{a^2_l|\dot
x_l(0)|^2}{(\pi n_l)^2}}.$$

If $a_{m1}=a_{m2}=\ldots=a_{mn}=a_m$, then the multiindex $\nb$
reduces to the index $k\in\mathbb N$. The equation~\eqref{an7}
implies
$$l_k^2=|\dot x(0)|^2=4\pi k\sum_{m=1}^{3}\frac{z_m^2(1)}{a_m^2},\qquad k\in\mathbb N.$$

Let $U$ be a neighborhood of the origin $O$. From Theorem~\ref{4},
we know that no matter how small $U$ is, we can always find points
in $U$ which are connected to $O$ by an infinite number of
geodesics. This is totally different from the Riemannian geometry.
It is known that every point of a Riemannian manifold is connected
to every other point in a sufficiently small neighborhood by one
single, unique geodesic.
\end{remark}

\subsection{Connectivity between $(0,0)$ and $(x,z)$, $x\neq 0$, $z\neq 0$.}

Now, we will look for a solution of the equation~\eqref{9} with
the boundary conditions
$$x(0)=0,\quad z(0)=0,\qquad x(1)=x,\quad z(1)=z.$$

Let us make some preliminary calculations. We obtain
\begin{equation}\label{66}
|\dot
x_l(0)|^2=\frac{|\theta|_l^2}{\sin^2|\theta|_l}|x_l(1)|^2,\quad
l=1,\ldots,n,
\end{equation} from~\eqref{48} for $s=1$.
Putting $s=1$ in~\eqref{53} and making use of~\eqref{66} we obtain
\begin{equation}\label{61}
z_m (1)=\frac{\theta_m}{4}\sum_{l=1}^{n}\frac{a_{ml}^2|
x_l(1)|^2}{|\theta|_l}\mu(|\theta|_l),\quad m=1,2,3.
\end{equation}
where
$\mu(|\theta|_l)=\frac{|\theta|_l}{\sin^2(|\theta|_l)}-\cot(|\theta|_l)$.
The function $\mu(\theta)$, introduced by Gaveau in~\cite{Ga}, was
first studied in detailed by Beals, Gaveau, Greiner
in~\cite{BGGr2,BGGr1,CChGr1}. By the following lemma, one finds
some basic properties of the function $\mu$.

\begin{figure}[ht]
\centering \scalebox{0.7}{\includegraphics{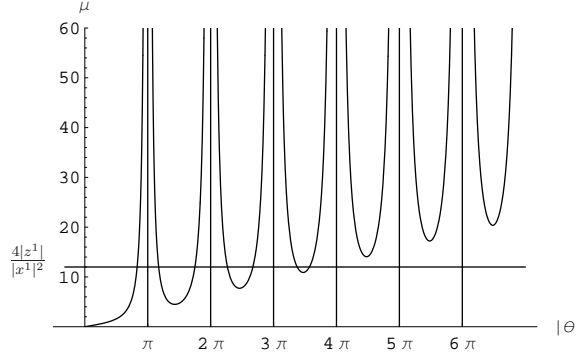}}
\caption[]{The graph of $\mu(\theta)$ and solutions of the
equation $\mu(|\theta|)=\frac{4|z(1)|}{|x(1)|^2}$}\label{fig1}
\end{figure}

\begin{lemma}\label{56}
The function $\mu(\theta)=\frac{\theta}{\sin^2\theta}-\cot\theta$
is an increasing diffeomorphism of the interval $(-\pi,\pi)$ onto
$\mathbb R$. On each interval $(m\pi,(m+1)\pi)$, $m=1,2,\ldots$,
the function $\mu$ has a unique critical point $c_m$. On this
interval the function $\mu$ strictly decreases from $+\infty$ to
$\mu(c_m)$, and then, strictly increases from $\mu(c_m)$ to
$+\infty$. Moreover,
$$\mu(c_m)+\pi<\mu(c_{m+1}),\qquad m=1,2,\ldots.$$
\end{lemma}
The graph of $\mu(\theta)$ is given in Figure~\ref{fig1}.
%%%%%%%%%%%%%%%%%%%%%%%%FIGURE 1%%%%%%%%%%%%%%%%%%%%%%%%%%%%%%
%%%%%%%%%%%%%%%%%%%%%%%%%%%%%%%%%%%%%%%%%%%%%%%%%%%%%%%%%%%%%%

\begin{theorem}\label{3}
Given a point $P(x,z)$ with $x_l\neq 0$, $l=1,\ldots,n$, $z\neq
0$, there are finitely many geodesics joining the point $O(0,0)$
with a point $P$. Let
$\vartheta_{(1)}=(|\theta_1|_1,\ldots,|\theta_1|_n),\ldots,\vartheta_{(N)}=(|\theta_N|_1,\ldots,|\theta_N|_n)$
be solutions of the system
\begin{equation}\label{67}\sum_{m=1}^{3}
\frac{16z_m^2(1)a_{ml}^2}{\Big(\sum_{r=1}^{n}\frac{a_{mr}^2|x_r(1)|^2\mu(|\theta|_r)}{|\theta|_r}\Big)^2}
=|\theta|_l^2,\quad l=1,\ldots,n.\end{equation} We fix one of the
solution $\vartheta=(|\theta|_1,\ldots,|\theta|_n)$. Then the
equations of the geodesics are
\begin{equation}\label{57}
\begin{split}
x^{(\nb)}_l(s) = & \big(4\cot(|\theta|_l)\sin^2(s|\theta|_l)
-2\sin(2s|\theta|_l)\big)\frac{[\mathbf Z]_l}{|\theta|_l}x_l(1)
\\ & +
\big(\frac{1}{2}\cot|\theta|_l\sin(2s|\theta|_l)+\sin^2(s|\theta|_l)\big)\U
x_l(1),\quad l=1,\ldots,n,\quad \nb=1,2,\ldots,N,
\end{split}
\end{equation}
\begin{equation}\label{58}
z^{(\nb)}_m(s) =z_m(1) \frac{\sum_{l=1}^{n}
\frac{a_{ml}^2|x(1)|^2_l}{\sin^2(|\theta|_l)}\big(s-\frac{\sin(2s|\theta|_l}{2|\theta|_l})\big)}{\sum_{l=1}^{n}\frac{a_{ml}^2|x(1)|^2_l}{|\theta|_l}\mu(|\theta|_l)},
\quad m=1,2,3,\quad \nb=1,2,\ldots,N,
\end{equation} where $\mathbf Z$ is a block diagonal matrix with
blocks~\eqref{an6}. The lengths of these geodesics are
\begin{equation}\label{an17}l^2_{\nb}=\sum_{l=1}^{n}\frac{|\theta|_l^2|x_l(1)|^2}{\sin^2(|\theta|_l)}
=16\sum_{m=1}^{3}\frac{z^2_m(1)}{\sum_{l=1}^{n}\frac{a_{ml}^2|x_l(1)|^2}{|\theta|_l}\mu(|\theta|_l)}+
\sum_{l=1}^{n}|x_l(1)|^2|\theta|_l\cot(|\theta|_l).\end{equation}
\end{theorem}
\begin{proof}
We have \begin{equation}\label{an5}\theta_m=\frac{4
z_m(1)}{\sum_{l=1}^{n}\frac{a_{ml}^2|
x_l(1)|^2}{|\theta|_l}\mu(|\theta|_l)},\quad
m=1,2,3,\end{equation} from~\eqref{61}. Then
\begin{equation}\label{an10}|\theta|^2_l=\sum_{m=1}^{3}\theta^2_m
a_{ml}^2=\sum_{m=1}^{3}
\frac{16z_m^2a_{ml}^2}{\Big(\sum_{r=1}^{n}\frac{a_{mr}^2|x_r(1)|^2\mu(|\theta|_r)}{|\theta|_r}\Big)^2},\quad
l=1,\ldots,n,\end{equation} that prove~\eqref{67}.

Let us fix one of the solutions of the equation~\eqref{67}
$\vartheta=(|\theta|_1,\ldots,|\theta|_n)$ for a given point
$P(x,z)$. Putting~\eqref{66} and~\eqref{an5} into~\eqref{53}, we
obtain~\eqref{58}.

Setting $s=1$ in~\eqref{48}, we find $\dot x^{(\nb)}_l(0)$ for
$\vartheta=(|\theta|_1,\ldots,|\theta|_n)$:
$$\dot
x^{(\nb)}_l(0)=2|\theta|_l\Big[\sin(2|\theta|_l)\U+\big(1-\cos(2|\theta|_l)\big)
\frac{[\Mb]_l}{|\theta|_l}\Big]^{-1}x_l(1)
=\Big[\big(|\theta|_l\cot|\theta|_l\big)\U-[\Mb]_l\Big]x_l(1).$$
This and~\eqref{48} imply
\begin{eqnarray}\label{70}
x^{(\nb)}_l(s)& = &
\frac{1}{2}\Big[\big(2\cot|\theta|_l\sin^2(s|\theta|_l)
-\sin(2s|\theta|_l)\big)\frac{[\Mb]_l}{|\theta|_l}\nonumber \\ & +
&
\big(\cot|\theta|_l\sin(2s|\theta|_l)+2\sin^2(s|\theta|_l)\big)\U\Big]x_l(1),\qquad
l=1,\ldots,n.
\end{eqnarray} Taking into account~\eqref{an5}, we deduce that each
block $[\Mb]_l$ of the matrix $\Mb$ takes the form $4[\mathbf Z]_l$
with block $[\mathbf Z]_l$ written as
\begin{equation}\label{an6} \left[\array{cccc} 0 &
\frac{z_{1}(1)a_{1l}}{\sum_{l=1}^{n}\frac{a_{1l}^2|
x_l|^2}{|\theta|_l}\mu(|\theta|_l)} &
-\frac{z_{3}(1)a_{3l}}{\sum_{l=1}^{n}\frac{a_{3l}^2|
x_l|^2}{|\theta|_l}\mu(|\theta|_l)} &
-\frac{z_{2}(1)a_{2l}}{\sum_{l=1}^{n}\frac{a_{2l}^2|
x_l|^2}{|\theta|_l}\mu(|\theta|_l)}
\\
-\frac{z_{1}(1)a_{1l}}{\sum_{l=1}^{n}\frac{a_{1l}^2|
x_l|^2}{|\theta|_l}\mu(|\theta|_l)} & 0 &
-\frac{z_{2}(1)a_{2l}}{\sum_{l=1}^{n}\frac{a_{2l}^2|
x_l|^2}{|\theta|_l}\mu(|\theta|_l)} &
\frac{z_{3}(1)a_{3l}}{\sum_{l=1}^{n}\frac{a_{3l}^2|
x_l|^2}{|\theta|_l}\mu(|\theta|_l)}
\\
\frac{z_{3}(1)a_{3l}}{\sum_{l=1}^{n}\frac{a_{3l}^2|
x_l|^2}{|\theta|_l}\mu(|\theta|_l)} &
\frac{z_{2}(1)a_{2l}}{\sum_{l=1}^{n}\frac{a_{2l}^2|
x_l|^2}{|\theta|_l}\mu(|\theta|_l)} & 0 &
\frac{z_{1}(1)a_{1l}}{\sum_{l=1}^{n}\frac{a_{1l}^2|
x_l|^2}{|\theta|_l}\mu(|\theta|_l)}
\\
\frac{z_{2}(1)a_{2l}}{\sum_{l=1}^{n}\frac{a_{2l}^2|
x_l|^2}{|\theta|_l}\mu(|\theta|_l)} &
-\frac{z_{3}(1)a_{3l}}{\sum_{l=1}^{n}\frac{a_{3l}^2|
x_l|^2}{|\theta|_l}\mu(|\theta|_l)} &
-\frac{z_{1}(1)a_{1l}}{\sum_{l=1}^{n}\frac{a_{1l}^2|
x_l|^2}{|\theta|_l}\mu(|\theta|_l)} & 0\endarray\right],
\end{equation}
Finally,~\eqref{an6} and~\eqref{70} give~\eqref{57}.

To obtain the length of the geodesics, we make the following
calculations.
\begin{eqnarray}\label{an8}
\sum_{m=1}^{3}z_m(1)\theta_m & = &
\frac{1}{4}\sum_{l=1}^{n}\frac{|x_l(1)|^2\mu(|\theta|_l)}{|\theta|_l}
\sum_{m=1}^{3}\theta^2_ma^2_{ml}=\frac{1}{4}\sum_{l=1}^{n}|x_l(1)|^2|\theta|_l\mu(|\theta|_l)\nonumber
\\ & = & \frac{1}{4}\sum_{l=1}^{n}\frac{|x_l(1)|^2|\theta|_l^2}{\sin^2(|\theta|_l)}
-\frac{1}{4}\sum_{l=1}^{n}|x_l(1)|^2|\theta|_l\cot(|\theta|_l)
\\ & = & \frac{|\dot
x(0)|^2}{4}-\frac{1}{4}\sum_{l=1}^{n}|x_l|^2|\theta|_l\cot(|\theta|_l).\nonumber
\end{eqnarray}
From the other hand~\eqref{an5} implies
\begin{equation}\label{an9}\sum_{m=1}^{3}z_m(1)\theta_m = 4
\sum_{m=1}^{3}\frac{z^2_m(1)}{\sum_{l=1}^{n}\frac{a_{ml}^2|x_l(1)|^2}{|\theta|_l}\mu(|\theta|_l)}.
\end{equation}
The formula~\eqref{an17} follows from~\eqref{an8} and~\eqref{an9}.
\end{proof}

\begin{remark}
Let us consider the particular case $a_{1l}=a_{2l}=a_{3l}=a_l>0$.
We have
\begin{equation}\label{an18}
4|z|=\sum_{l=1}^{n}|a_l||x_l(1)|^2\mu(a_l|\theta|)\end{equation}
from~\eqref{an10}. Here
$|\theta|^2=\theta_1^2+\theta_2^2+\theta_3^2$. We denote by
$|\theta|_1,\ldots,|\theta|_N$ the solutions of~\eqref{an18} and
let $|\theta|$ one of the solutions. Then~\eqref{an5} implies
$$\theta_m=\frac{4z_m(1)|\theta|}{\sum_{l=1}^{n}a_l|x_l|^2\mu(a_l|\theta|)}=\frac{z_m}{|z|}|\theta|.$$
To obtain the formula for the length of geodesics, we
write~\eqref{an18} as
$$4|z|=\frac{1}{|\theta|}\sum_{l=1}^{n}\frac{|x_l|^2|\theta|_l^2}{\sin^2(|\theta|_l)}
-\sum_{l=1}^{n}|a_l||x_l|^2\cot(|\theta|_l)=\frac{l^2_{\nb}}{|\theta|}
-\sum_{l=1}^{n}|a_l||x_l|^2\cot(|\theta|_l)$$ and get
$$l^2_{\nb}=|\theta|\big(4|z(1)|+\sum_{l=1}^{n}a_l|x_l(1)|^2\cot(a_l|\theta|)\big).$$

Simplifying more the situation and supposing that $a_l=a>0$ for
all $l=1,\ldots,n$, we get that $|\theta|_l=a|\theta|$. This
implies that $|\theta|$ is a solution of the equation (see
Figure~\ref{fig1})
\begin{equation}\label{an11}\mu(a|\theta|)=\frac{4|z(1)|}{a|x(1)|^2}.\end{equation} In this case to calculate
the length of geodesics joining $(0,0)$ to $(x,z)$, $x_l\neq 0$,
$l=1,\ldots,n$, we use the homogeneous norm
$|(x,z)|^2=|x|^2+4|z|$. It gives for a solution
$|\theta|_{\alpha}$, $\alpha=1,\ldots,N$, of~\eqref{an11}
$$|x|^2+4|z|=|x|^2+a|x|^2\mu(a|\theta|_{\alpha})=(1+a\mu(a|\theta|_{\alpha}))
\frac{\sin^2(a|\theta|_{\alpha})}{a^2|\theta|_{\alpha}^2}l^2_{\alpha}$$
by~\eqref{an11} and~\eqref{66}. Then
$$l^2_{\alpha}=
\frac{a^2|\theta|_{\alpha}^2}{\sin(a|\theta|_{\alpha})\big(\sin(a|\theta|_{\alpha})-\cos(a|\theta|_{\alpha})\big)+a^2|\theta|_{\alpha}}
(|x|^2+4|z|).$$

In the last simplest case it is easy to observe that if $z$ is
fixed, and $|x|$ tends to zero, then the ratio
$\frac{4|z|}{a|x|^2}$ increases and the number of solutions of the
equation $\frac{4|z|}{a|x|^2}=\mu(a|\theta|)$ also increases (see
Figure~\ref{fig1}). In this case, the function
$\mu(a|\theta|_{\alpha})=\frac{a|\theta|_{\alpha}-\cos(a|\theta|_{\alpha})\sin(a|\theta|_{\alpha})}{\sin^2(a|\theta|_{\alpha})}$
tends to infinity as $|x|\to 0$, and we obtain that
$\sin^2(a|\theta|)=0$ and $a|\theta|=\pi n$, $n\in\mathbb N$. One
sees that Theorem~\ref{4} is the limiting case of Theorem~\ref{3}
as the ratio $\frac{4|z|}{a|x|^2}$ tends to $\infty$. If we fix
$x$ and let $|z|$ tend to $0$, then the equation
$\frac{4|z|}{a|x|^2}=\mu(a|\theta|)$ says that $\mu(a|\theta|)\to
0$. This implies that $|\theta|\to 0$ and we obtain
Theorem~\ref{54} as another limit case of Theorem~~\ref{3}.

The last particular case is when
$a_{m1}=a_{m2}=\ldots=a_{mn}=a_m$. We denote
$|\theta|_l^2=a_1^2\theta_1^2+a_2^2\theta_2^2+a_3^2\theta_3^2=|\theta|_a^2$.
The equation to find $|\theta|_a$ is
$$\mu(|\theta|_a)=\frac{4}{|x(1)|^2}\sqrt{\frac{z_m^2(1)}{a_m^2}}.$$
The value of $\theta_m$ and the lengths of geodesics are
$$\theta_m=\frac{z_m(1)|\theta|_a}{a^2_m\sqrt{\frac{z_m^2(1)}{a_m^2}}},\qquad
l^2(|\theta|_a)=|\theta|_a\Big(4\sqrt{\frac{z_m^2(1)}{a_m^2}}+|x(1)|^2\cot(|\theta|_a)\Big).$$
\end{remark}

In the following theorem we consider the connection between the
origin and a point $P(x,z)$ when some of the coordinates $x_l$
vanish.

\begin{theorem}\label{an19}
Given a point $P(x,z)$ with $x_l\neq 0$, $l=1,\ldots,p-1$, and
$x_l=0$, $l=p,\ldots,n$, $z\neq 0$, there are infinitely many
geodesics joining the point $O(0,0)$ with a point $P$. Let
$S_{1m}=\sum_{l=1}^{p-1}\frac{a_{ml}^2|x_l(1)|^2}{|\theta|_l}\mu(|\theta|_l)$,
$S_{2m}=\sum_{l=p}^{n}\frac{a_{ml}^2|\dot x_l(0)|^2}{\pi^2
n^2_l}$, $\nb_{\beta}=(n_p,\ldots,n_n)$ be a multiindex with
positive integer-valued components for each $\beta\in\mathbb N$,
and $\vartheta_{\kappa}=(|\theta|_1,\ldots,|\theta|_{p-1})$,
$\kappa=1,\ldots,N$ be solutions of the system
\begin{equation}\label{an23}|\theta|_l^2=\sum_{m=1}^{3}\frac{16z_m^2(1)a_{ml}^2}{\big(S_{1m}+S_{2m}\big)^2},\quad
l=1,\ldots,p-1.\end{equation} Then the equations of geodesics are
\begin{equation}\label{an24}
\begin{split}
x^{(\kappa)}_l(s) = & \big(4\cot(|\theta|_l)\sin^2(s|\theta|_l)
-2\sin(2s|\theta|_l)\big)\frac{[\mathbf Z]_l}{|\theta|_l}x_l(1)
\\ & +
\big(\frac{1}{2}\cot|\theta|_l\sin(2s|\theta|_l)+\sin^2(s|\theta|_l)\big)\U
x_l(1),\quad l=1,\ldots,n,\quad \kappa=1,2,\ldots,N,
\end{split}\end{equation}
\begin{equation}\label{an241} x^{(\nb_{\beta})}_l(s) = 2\frac{1-\cos(2s\pi n_l)}{(\pi
n_l)^2}[\mathbf Z]_l\dot x_l(0) +\frac{\sin(2s\pi n_l)}{2\pi
n_l}\U\dot x_l(0),\ \ l=1,\ldots,n,\quad\beta\in\mathbb
N,\end{equation} where
\begin{equation}\label{an26} [\mathbf Z]_l= \left[\array{cccc} 0 &
\frac{z_{1}a_{1l}}{S_{11}+S_{21}} &
-\frac{z_{3}a_{3l}}{S_{13}+S_{23}} &
-\frac{z_{2}a_{2l}}{S_{12}+S_{22}}
\\
-\frac{z_{1}a_{1l}}{S_{11}+S_{21}} & 0 &
-\frac{z_{2}a_{2l}}{S_{12}+S_{22}} &
\frac{z_{3}a_{3l}}{S_{13}+S_{23}}
\\
\frac{z_{3}a_{3l}}{S_{13}+S_{23}} &
\frac{z_{2}a_{2l}}{S_{12}+S_{22}} & 0 &
\frac{z_{1}a_{1l}}{S_{11}+S_{21}}
\\
\frac{z_{2}a_{2l}}{S_{12}+S_{22}} &
-\frac{z_{3}a_{3l}}{S_{13}+S_{23}} &
-\frac{z_{1}a_{1l}}{S_{11}+S_{21}} & 0\endarray\right],
\end{equation}
and
\begin{equation}\label{an25}
z^{(\kappa,\nb_{\beta})}_m(s) =\frac{z_m(1)}{S_{1m}+S_{2m}}
\Big(\sum_{l=1}^{p-1}
\frac{a_{ml}^2|x(1)|^2_l}{\sin^2(|\theta|_l)}\big(s-\frac{\sin(2s|\theta|_l)}{2|\theta|_l})\big)+
\sum_{l=p}^{n}\frac{a_{ml}^2|\dot x_l(0)|^2}{\pi^2
n_l^2}\big(s-\frac{\sin(2s\pi n_l)}{2\pi n_l})\big)\Big),
\end{equation} with $m=1,2,3$.

The lengths of these geodesics are
\begin{equation}\label{an27}l^2_{\kappa,\nb_{\beta}}=\sum_{l=1}^{n}|\dot x_l(0)|^2
=16\sum_{m=1}^{3}\frac{z^2_m(1)}{S_{1m}+S_{2m}}+
\sum_{l=1}^{p-1}|x_l(1)|^2|\theta|_l\cot(|\theta|_l).\end{equation}
\end{theorem}
\begin{proof} If $x_l(1)=0$, then the formula $|x_l(1)|^2=\frac{\sin^2(|\theta|_l)}{|\theta|_l^2}|\dot
x_l(0)|^2$ implies that $|\dot x_l(0)|=0$ or
$\sin^2(|\theta|_l)=0$. If $|\dot x_l(0)|=0$, then the
corresponding $x_l(s)\equiv 0$. The more interesting case when
$|\dot x_l(0)|\neq 0$ for $l=p,\ldots,n$. Then $|\theta|_l=\pi
n_l$, $n_l\in\mathbb N$, $l=p\ldots,n$. We deduce from~\eqref{53}
for $s=1$
\begin{equation}\label{an20}
z_m(1)=\frac{\theta_m}{4}\Big(\sum_{l=1}^{p-1}\frac{a_{ml}^2|x_l(1)|^2}{|\theta|_l}\mu(|\theta|_l)
+\sum_{l=p}^{n}\frac{a_{ml}^2|\dot x_l(0)|^2}{\pi^2
n^2_l}\Big)=\frac{\theta_m}{4}(S_{1m}+S_{2m}),
\end{equation} where the number $n_l$ can have any positive
integer value. We conclude that the sum
$S_{2m}=\sum_{l=p}^{n}\frac{a_{ml}^2|\dot x_l(0)|^2}{\pi^2 n^2_l}$
admits countably many values. To define $|\theta|_l$ we find
$\theta_m=\frac{4z_m(1)}{S_{1m}+S_{2m}}$, $m=1,2,3$,
from~\eqref{an20} and then argue as in Theorem~\ref{3}:
\begin{equation}\label{an21}
|\theta|_l^2=\sum_{m=1}^{3}\theta_m^2a^2_{ml}=\sum_{m=1}^{3}\frac{16z_m^2(1)a_{ml}^2}{\big(S_{1m}+S_{2m}\big)^2},\quad
l=1,\ldots,p-1.\end{equation} Conclude, that for each multiindex
with positive integer-valued components
$\nb_{\beta}=(n_{p},\ldots,n_n)$, $\beta\in\mathbb N$, the
equation~\eqref{an21} defines the multiindex
$\vartheta_{\kappa}=(|\theta|_1,\ldots,|\theta|_{p-1})$,
$\kappa=1,\ldots, N$. Let us fix one of the solutions
$(\vartheta_{\kappa},\nb_{\beta})=(|\theta|_1,\ldots,|\theta|_{p-1},n_{p},\ldots,n_n)$.
The relations~\eqref{48} for $s=1$ give
\begin{equation}\label{an22}
\dot x_l(0)=\Big((|\theta|_l\cot
|\theta|_l)\U-[\Mb]_l\Big)x_l(1),\quad l=1,\ldots,p-1.
\end{equation}
Substituting~\eqref{an22} into~\eqref{48} and~\eqref{53}, we
obtain~\eqref{an24},~\eqref{an241}, and~\eqref{an25}. We get the
form of~\eqref{an26} from
$\theta_m=\frac{4z_m(1)}{S_{1m}+S_{2m}}$, $m=1,2,3$ and the
definition of the matrix $\Mb$.

To calculate the length of the geodesic we argue as follows:
\begin{eqnarray}\label{an28}
\sum_{m=1}^{3}z_m\theta_m & = &
\frac{1}{4}\Big(\sum_{m=1}^{3}\theta^2_ma^2_{ml}\sum_{l=1}^{p-1}\frac{|x_l(1)|^2\mu(|\theta|_l)}{|\theta|_l}
+\sum_{m=1}^{3}\theta^2_ma^2_{ml}\sum_{l=p}^{n}\frac{|\dot
x_l(0)|^2}{\pi^2 n_l^2}\Big)\nonumber
\\ & = & \frac{1}{4}\sum_{l=1}^{p-1}|x_l(1)|^2|\theta|_l\mu(|\theta|_l)+\sum_{l=p}^{n}\frac{|\dot
x_l(0)|^2}{4}\nonumber
\\ & = & \frac{1}{4}\sum_{l=1}^{p-1}\frac{|x_l(1)|^2|\theta|_l^2}{\sin^2(|\theta|_l)}
-\frac{1}{4}\sum_{l=1}^{p-1}|x_l(1)|^2|\theta|_l\cot(|\theta|_l)+\sum_{l=p}^{n}\frac{|\dot
x_l(0)|^2}{4}
\\ & = & \frac{|\dot
x(0)|^2}{4}-\frac{1}{4}\sum_{l=1}^{p-1}|x_l(1)|^2|\theta|_l\cot(|\theta|_l).\nonumber
\end{eqnarray}
From the other hand, since
$\theta_m=\frac{4z_m(1)}{S_{1m}+S_{2m}}$, $m=1,2,3$,  we deduce
\begin{equation}\label{an29}\sum_{m=1}^{3}z_m\theta_m = 4
\sum_{m=1}^{3}\frac{z^2_m(1)}{S_{1m}+S_{2m}}.
\end{equation}
The formula~\eqref{an27} follows from~\eqref{an28}
and~\eqref{an29}.
\end{proof}

\begin{remark} Let make some simulations for the anisotropic group $Q^2$.
Set $$x_1(1)=(x_{11},x_{12},0,0),\quad x_2(1)=0,\quad \dot
x_2(0)=(\dot x_{21}(0),\dot x_{22}(0),0,0),$$ $$z_1(1)\neq
0,\qquad z_2(1)=z_3(1)=0.$$ In this case the equation~\eqref{an23}
can be written in the form
\begin{equation}\label{an31}\mu(|\theta|_1)=\frac{4|z_1|}{a_{11}|x_1(1)|^2}-\frac{|\theta|_1a^2_{12}|\dot
x_2(0)|^2}{\pi^2 n^2 a^2_{11}|x_1(1)|^2}.\end{equation} We present
the solutions for different values of $n$: $n=1,2,50$ in
Figure~\ref{an30}. We see that for sufficiently big value of $n$
the second term in the right hand side of~\eqref{an31} goes to
$0$, and we obtain a finite number of solutions for $|\theta|_1$.

%%%%%%%%%%%%%%%%%%%%%%%%FIGURE 1%%%%%%%%%%%%%%%%%%%%%%%%%%%%%%
\begin{figure}[ht]
\centering \scalebox{0.6}{\includegraphics{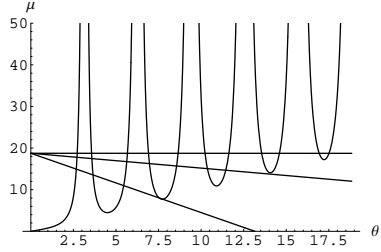}}
\caption[]{Solutions of equation~\eqref{an31} \label{an30}}
\end{figure}
%%%%%%%%%%%%%%%%%%%%%%%%%%%%%%%%%%%%%%%%%%%%%%%%%%%%%%%%%%%%%%

Moreover, we obtain countably many geodesics, because of the
second part of multiindex, corresponding to the positive integer
values is countably infinite.  Nevertheless, since the sums
$S_{2m}$ tends to $0$ and the sums $S_{1m}$ are strictly positive
as $\nb\to \infty$, we conclude that the lengths of these
geodesics are bounded from the above. The projection in each
subspace $x_l$ are still ellipsoids. In Figures~\ref{an32}
and~\ref{an33} we present the projection of a geodesic into spaces
$(x_{1},z_1)$ and $(x_2,z_1)$. We can see that the number of loops
is different and increases in the subspace corresponding vanishing
value of $x_{l}(1)$.

%%%%%%%%%%%%%%%%%%%%%%%%FIGURE 1%%%%%%%%%%%%%%%%%%%%%%%%%%%%%%
\begin{figure}[ht]
\centering \scalebox{0.55}{\includegraphics{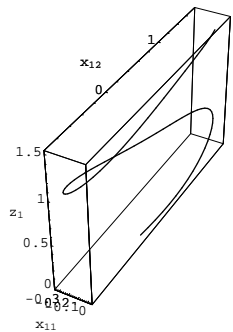}}
\caption[]{Projection of a geodesic to the space
$(x_{1},z_1)$\label{an32}}
\end{figure}
%%%%%%%%%%%%%%%%%%%%%%%%%%%%%%%%%%%%%%%%%%%%%%%%%%%%%%%%%%%%%%

%%%%%%%%%%%%%%%%%%%%%%%%FIGURE 1%%%%%%%%%%%%%%%%%%%%%%%%%%%%%%
\begin{figure}[ht]
\centering \scalebox{0.55}{\includegraphics{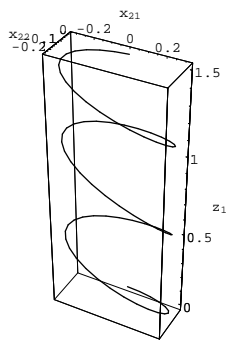}}
\caption[]{Projection of a geodesic to the space
$(x_{2},z_1)$\label{an33}}
\end{figure}
%%%%%%%%%%%%%%%%%%%%%%%%%%%%%%%%%%%%%%%%%%%%%%%%%%%%%%%%%%%%%%
\end{remark}

\section{Complex Hamiltonian mechanics}

Our aim now is to study the complex action which may  be used to obtain
the length of real geodesics.

\begin{definition}
A complex geodesic is the projection  of a solution of the
Hamiltonian system~\eqref{2} with the non-standard boundary
conditions $$x(0)=0,\quad x(1)=x,\quad z(0)=0,\quad
z(1)=z,\quad\text{and}$$
$$\theta_{m}=-i\tau_m,\quad m=1,2,3,$$
on the $(x,z)$-space.
\end{definition}
Let us introduce the notation $-i\tau$ for the vector
$(-i\tau_1,-i\tau_2,-i\tau_3)$. We write
$|\tau|_l=\sqrt{a_{1l}^2\tau_1^2+a_{2l}^{2}\tau_2^2+a_{3l}^{3}\tau_3^2}$.
Then
$|\theta|_l=\sqrt{a_{1l}^{2}\theta_{1}^2+a_{2l}^{2}\theta_{2}^2+a_{3l}^{2}\theta_{3}^2}=i|\tau|_l$.

Notice, that we should treat the missing directions apart from the
directions in the underlying space.

\begin{definition}
The modifying complex action is defined as
\begin{equation}
f(x,z,\tau)=-i\sum_m\tau_mz_{m}+\int_{0}^{1}\big((\dot
x,\xi)-H(x,z,\xi,\tau)\big)\,ds.
\end{equation}
\end{definition}

We present some useful calculations following from the
system~\eqref{2}.
\begin{equation}\label{36}
\begin{split}
(\xi,\dot x) &= 2|\xi|^2+(\Mb x,\xi)= \frac{1}{2}|\dot
x|^2-\frac{1}{2}(\Mb x,\dot x),
\\ |\xi|^2 = & \frac{|\dot x|^2}{4}-\frac{1}{2}(\Mb x,\dot
x)+\frac{1}{4}(\Mb x,\Mb x) =  \frac{|\dot
x|^2}{4}-\frac{1}{2}(\Mb x,\dot x)+\frac{1}{4}(\Tb^2 x,x),
\\ (\Mb x,\xi) = & \frac{1}{2}(\Mb x,\dot
x)-\frac{1}{2}(\Tb^2 x,x).
\end{split}
\end{equation}

Making use of the formulas~\eqref{35}, ~\eqref{36},
and~\eqref{66}, we deduce
\begin{equation*}
\begin{split}
f(x,z,\tau) = & -i\sum_m\tau_mz_{m} +\int_{0}^{1}\big((\dot
x,\xi)-H(x,z,\xi,\tau)\big)\,ds
\\ = & -i\sum_m\tau_mz_{m}+\int_{0}^{1}\Big(\frac{|\dot x(s)|^2}{4}-\frac{1}{2}(\Mb x,\dot
x)\Big)\,ds \\ = & -i\sum_m\tau_mz_{m}+\sum_{l=1}^{n}\frac{|\dot
x_l(0)|^2}{4}\int_{0}^{1}\cosh(2s|\tau|_l)\,ds
\\ =&-i\sum_m\tau_mz_{m}+\sum_{l=1}^{n}\frac{|x_l|^2}{4}\frac{(i|\tau|_l)^2}{\sin^2(-i|\tau|_l)}
\frac{\sinh(2|\tau|_l)}{2|\tau|_l}\\
=&-i\sum_m\tau_mz_{m}+\sum_{l=1}^{n}\frac{|x_l|^2}{4}|\tau|_l\coth|\tau|_l.
\end{split}
\end{equation*}

The complex action function satisfies the Hamilton-Jacobi equation
\begin{equation}
\label{eq:HJ} \sum_{m=1}^{3}\tau_m\frac{\partial f}{\partial
\tau_m} +H(x,z,\nabla_x f,\nabla_z f)=f.
\end{equation}
Indeed, we have
$$\frac{\partial f}{\partial
\tau_m}=-iz_{m}-i\tau_m\sum_{l=1}^{n}\frac{a_{ml}^2|x_l|^2}{4|\tau|_l}\mu(i|\tau|_l),\quad\
m=1,2,3.$$
$$H(x,z,\frac{\partial f}{\partial x},\frac{\partial f}{\partial
z})=H(x,z,\xi,\tau)=\sum_{l=1}^{n}\frac{|x_l|^2}{4}\frac{|\tau|_l^2}{\sinh^2|\tau|_l}$$
from~\eqref{35}, ~\eqref{36}, and~\eqref{47}. Then,
\begin{equation*}
\begin{split}
\sum_{m=1}^{3}\tau_m\frac{\partial f}{\partial \tau_m} & +
H(x,z,\frac{\partial f}{\partial x},\frac{\partial f}{\partial
z})=-i\sum_m\tau_mz_m+\sum_{l=1}^{n}\frac{|x_l|^2|\tau|_l}{4}\Big(-i\mu(i|\tau|_l)
+\frac{|\tau|_l}{\sinh^2|\tau|_l}\Big)
\\ = &
-i\sum_m\tau_mz_m+\sum_{l=1}^{n}\frac{|x_l|^2}{4}|\tau|_l\coth|\tau|_l=f.
\end{split}
\end{equation*} In the critical points $\tau_c$, where
$\frac{\partial f}{\partial \tau_m}=0$ we have from~\eqref{an8}
$$f(x,z,\tau_c)=H(x,z,\nabla_xf,\nabla_zf)=\frac{\mathcal
E}{2}=\frac{l^2}{4}(\gamma),$$ where a geodesic curve $\gamma$
connects the origin with $(x,z)$.

\section{Green's function for the Schr\"{o}dinger operator}

Consider the Schr\"{o}dinger operator
$$L=\Delta_0-i\frac{\partial}{\partial u}.$$ We are looking for a
distribution $P=P(x,z,u)$ on $\mathbb R^{4n}_x\times \mathbb
R^3_z\times\mathbb R^{+}_u$ satisfying the following conditions
\begin{itemize}\label{h3}
\item[1)]{$LP=\Delta_0P-i\frac{\partial P}{\partial u}=0$ for
$u>0$,}\item[2)]{$\lim\limits_{u\to
0^+}P(x,z,u)=\delta(x)\delta(z)$},
\end{itemize} where $\delta$ stands for the Dirac distribution.

The next propositions are easily verified.
\begin{proposition}\label{h4}
For any smooth function $\varphi$ and any smooth vector fields
$X_1,\ldots,X_n$ we have $$\Delta
e^{\varphi}=e^{\varphi}(\Delta\varphi+|\nabla\varphi|^2),$$ where
$\Delta=\sum_{j=1}^nX^2_j$ and
$|\nabla\varphi|^2=\sum_{j=1}^2(X_j\varphi)^2$.
\end{proposition}

We recall that $X$ denotes the horizontal gradient
$X_{11},\ldots,X_{4\,n}$.
\begin{proposition}\label{h5}
$$H(x,z,\nabla_xf,\nabla_zf)=|Xf|^2
=\sum_{k=1}^4\sum_{l=1}^n\big(X_{kl}f\big)^2.$$
\end{proposition}

\begin{proposition}\label{prop:trans}
Let $V$ and $f$ be smooth functions of $x$, $z$, $\tau$, and
$X_1,\ldots,X_n$ smooth vector fields. Then for any number $p$, we
have the following identity:
\begin{equation}\label{an12}
\Delta(Vf^{-p})=(\Delta V)f^{-p}-pf^{-p-1}\Big[(\Delta f)V
+2(\nabla f)(\nabla V)\Big]+(-p)(-p-1)f^{-p-2}V|\nabla f|^2.
\end{equation}
\end{proposition}
\begin{proof} The formula~\eqref{an12} is obtained by the direct
calculation.
\end{proof}
Before we go further, let us make some calculations. We apply
Proposition~\ref{h4} to $e^{-\frac{if}{u}}$ and system~\eqref{62}
of horizontal vector fields $X_{kl}$, $k=1,\ldots,4$,
$l=1,\ldots,n$. Introducing the notation $\varphi=-\frac{if}{u}$,
we get $\Delta_0\varphi=-\frac{i}{u}\Delta_0f$,
$|X\varphi|^2=-\frac{1}{u^2}|X f|^2$, and
\begin{equation}\label{h6}
\Delta_0e^{\varphi}=
e^{\varphi}\Big(-\frac{i}{u}\Delta_0f-\frac{1}{u^2}|X f|^2\Big) =
\frac{e^{\varphi}V(\tau)u^{2n+3}}{u^{2n+4}V(\tau)}\Big(-i\Delta_0f-\frac{1}{u}|Xf|^2\Big).
\end{equation}
Hamilton-Jacobi equation~\eqref{eq:HJ}, Proposition~\ref{h5}, and
the equality~\eqref{h6} imply
\begin{equation}\label{h7}\Delta_0\frac{e^{\varphi}V(\tau)}{u^{2n+3}}=
\frac{e^{\varphi}V(\tau)}{u^{2n+4}}\Big(-i\Delta_0f-\frac{f}{u}+\frac{1}{u}\sum_{m=1}^{3}\tau_m\frac{\partial
f}{\partial \tau_m}\Big).\end{equation} Differentiating
$\frac{e^{\varphi}V(\tau)}{u^{2n+3}}$ with respect to $u$, we
obtain
\begin{equation}\label{h8}
-i\frac{\partial}{\partial
u}\Big(\frac{e^{\varphi}V(\tau)}{u^{2n+3}}\Big)=\frac{e^{\varphi}V(\tau)}{u^{2n+4}}\Big(\frac{f}{u}+i(2n+3)\Big).
\end{equation}
Summing~\eqref{h7} and~\eqref{h8}, we have
\begin{equation}\label{h9}
\Big(\Delta_0-i\frac{\partial}{\partial
u}\Big)\frac{e^{\varphi}V(\tau)}{u^{2n+3}}=i\frac{e^{\varphi}V(\tau)}{u^{2n+4}}
\Big((2n+3)-\Delta_0f-\frac{i}{u}\sum_{m=1}^{3}\tau_m\frac{\partial
f}{\partial \tau_m}\Big).
\end{equation}
We express $-\frac{i}{u}\sum_{m=1}^{3}\tau_m\frac{\partial
f}{\partial \tau_m}$ from the formula
$$\sum_{m=1}^{3}\frac{\partial }{\partial \tau_m}\big(e^{\varphi}V(\tau)\tau_m\big)=
e^{\varphi}V(\tau)\Big(-\frac{i}{u}\sum_{m=1}^{3}\tau_m\frac{\partial
f}{\partial
\tau_m}\Big)+e^{\varphi}\sum_{m=1}^{3}\tau_m\frac{\partial
V(\tau)}{\partial \tau_m}+3e^{\varphi}V(\tau)$$ and put it
into~\eqref{h9}. Finally, we deduce
\begin{equation}\label{h11}
\begin{split}
\Big(\Delta_0-i\frac{\partial}{\partial
u}\Big)\frac{e^{\varphi}V(\tau)}{u^{2n+3}}  = &
i\frac{e^{\varphi}}{u^{2n+4}} \Big((2n-\Delta
f)V-\sum_{m=1}^{3}\tau_m\frac{\partial V}{\partial\tau_m}\Big)
\\  + & \frac{i}{u^{2n+4}}\sum_{m=1}^{3}\frac{\partial }{\partial
\tau_m}\big(e^{\varphi}V(\tau)\tau_m\big).\end{split}\end{equation}
The equation
\begin{equation}\label{h10}
(2n-\Delta f)V-\sum_{m=1}^{3}\tau_m\frac{\partial
V}{\partial\tau_m}=0\end{equation} is called the {\it transport
equation}. We show that the function
$$V(\tau)=\prod_{l=1}^{n}\frac{|\tau|_l^{2}}{\sinh^{2}|\tau|_l}$$ is a solution of
transport equation. Indeed, since
$$f=f(x,z,\tau)=-i\sum_m\tau_mz_{m}+\sum_{l=1}^{n}\frac{|x_l|^2}{4}|\tau|_l\coth(|\tau|_l),$$
we have
$$\frac{\partial f}{\partial z_m}=-i\tau_m,\quad \frac{\partial^2 f}
{\partial z^2_m}=0,\quad m=1,2,3.$$
$$\frac{\partial f}{\partial x_{kl}}=\frac{1}{2}x_{kl}|\tau|_l\coth(|\tau|_l),\quad
\frac{\partial^2 f}{\partial
x^2_{kl}}=\frac{|\tau|_l}{2}\coth(|\tau|_l),\quad k=1,\ldots,4,\
l=1,\ldots,n.$$ Finally,
$$\Delta f=2\sum_{l=1}^{n}|\tau|_l\coth(|\tau|_l)$$ and
\begin{equation}\label{h2}(2n-\Delta f)\,V(\tau)=
2V(\tau)\big(n-\sum_{l=1}^{n}|\tau|_l\coth(|\tau|_l)\big).\end{equation}

On the other hand the equalities \begin{eqnarray*}\frac{\partial
V}{\partial \tau_m} & =
&\sum\limits_{r=1}^{n}\prod\limits_{l=1,\,l\neq
r}^{n}\frac{|\tau|_l^{2}}{\sinh^{2}(|\tau|_l)}\cdot\frac{\partial
}{\partial \tau_m}\Big(\frac{|\tau|_r^{2}}{\sinh^{2}(|\tau|_r)}\Big)
\\ & = &\sum\limits_{r=1}^{n}\prod\limits_{l=1,\,l\neq
r}^{n}\frac{|\tau|_l^{2}}{\sinh^{2}(|\tau|_l)}\cdot\frac{2a_{mr}^2\tau_m}{\sinh^2|\tau|_r}\big(1-|\tau|_r\coth(|\tau|_r)\big)
,\quad m=1,2,3,\end{eqnarray*} imply
\begin{eqnarray*}\sum_{m=1}^{3}\tau_m\frac{\partial
V}{\partial\tau_m} & = &
\sum\limits_{r=1}^{n}\prod\limits_{l=1,\,l\neq
r}^{n}\frac{|\tau|_l^{2}}{\sinh^{2}(|\tau|_l)}\cdot\frac{2|\tau|^{2}_r}{\sinh^{2}(|\tau|_r)}
(1-|\tau|_r\coth(|\tau|_r))\\
& = &
2\prod\limits_{l=1}^{n}\frac{|\tau|_l^{2}}{\sinh^{2}(|\tau|_l)}
\cdot\sum\limits_{r=1}^{n}(1-|\tau|_r\coth(|\tau|_r))=
2V(\tau)(n-\sum\limits_{r=1}^{n}|\tau|_r\coth(|\tau|_r)),\end{eqnarray*}
that shows that $V(\tau)$ is a solution of the transport
equation~\eqref{h10}. The function $V(\tau)$ is called the {\it
volume element}.

If the volume element $V(\tau)$ satisfies the equation~\eqref{h10}
then the equation~\eqref{h11} is reduced to the next one
\begin{equation}\label{h12}
\Big(\Delta_0-i\frac{\partial}{\partial
u}\Big)\frac{e^{\varphi}V(\tau)}{u^{2n+3}}=
\frac{i}{u^{2n+4}}\sum_{m=1}^{3}\frac{\partial }{\partial
\tau_m}\big(e^{\varphi}V(\tau)\tau_m\big).\end{equation} We note
that the expression $e^{\varphi}V(\tau)\tau_m$ vanishes as
$|\tau|\to\infty$. Integrating over $\mathbb R^3$ with respect to
$d\tau=d\tau_1d\tau_2d\tau_3$, we obtain
$$\Big(\Delta_0-i\frac{\partial}{\partial
u}\Big)\int_{\mathbb
R^3}\frac{e^{\varphi}V(\tau)}{u^{2n+3}}\,d\tau=0\quad\text{for}\quad
u>0.$$ Thus, the function
$$P(x,z,u)=\frac{C}{u^{2n+3}}\int_{\mathbb
R^3}e^{\frac{-if}{u}}V(\tau)\,d\tau$$ satisfies the first
condition to the Green function at the origin of the
Schr\"{o}dinger operator.

\subsection{The heat kernel on $Q^n$}

In this section we denote the time variable by $t$ and we will
consider the heat operator $$\Delta_0-\frac{\partial}{\partial
t}=\sum_{k,l}Y_{kl}^2-\frac{\partial}{\partial t},$$ where
$Y=(Y_{1\,1},\ldots,Y_{4\,n})=\nabla_{y}+\frac{1}{2}(\sum_{m=1}^3\Mb_m
y\frac{\partial}{\partial w_m})$ with
$y=(y_{1\,1},\ldots,y_{4\,n})$. The fundamental solution at the
origin is the function $P(y,w,t)$ defined on $Q^n\times\mathbb
R^1_+$ such that the following conditions
\begin{itemize}
\item[1)]{$\Delta_0P-i\frac{\partial P}{\partial t}=0$\ \ \ for\ \
\ $t>0$,}\item[2)]{$\lim\limits_{t\to
0^+}P(y,w,t)=\delta(y)\delta(w)$}
\end{itemize} hold.
With the change of variables $$u=it,\quad x_{kl}=iy_{kl},\quad
w_m=z_m,\quad m=1,2,3,\ \ k=1,\ldots,4,\ \ l=1,\ldots,n,$$ the
heat operator transforms to the Schr\"{o}dinger operator
$$\Delta_0-i\frac{\partial}{\partial
u}=\sum_{k,l}X_{kl}^2-i\frac{\partial}{\partial u}.$$ Indeed,
under this change of variables we obtain
$$-i\frac{\partial}{\partial u}=\frac{\partial}{\partial
t}\qquad\text{and}\qquad X=-iY.$$

The calculus of the previous subsection give us the following
statement.
\begin{theorem}
\label{th:heat1}
The heat kernel at the origin is given by
$$P(y,w,t)=\frac{C}{t^{2n+3}}\int_{\mathbb
R^3}e^{\frac{-f}{t}}V(\tau)\,d\tau,$$ where $$f(y,w,\tau) =
-i\sum_m\tau_mw_{m}+\sum_{l=1}^{n}\frac{|y_l|^2}{4}|\tau|_l\coth(|\tau|_l)$$
is the modified complex action and
$$V(\tau)=\prod_{l=1}^{n}\frac{|\tau|_l^{2}}{\sinh^{2}(|\tau|_l)}$$ is the volume
element.
\end{theorem}

\subsection{Green function for the sub-Laplace operator.}
Let us integrate the kernel $P(x,z,u)$ with respect to the time
variable $u$ on $(0,\infty)$. That is
\begin{equation*}
\begin{split}
\int_0^{\infty} P(x,z,u)\,du= & \int_0^{\infty} \frac{C}{u^{2n+3}}
\int_{\mathbb R^3} e^{-if/u}V(\tau)\, d\tau\ du\\ = &
C\int_{\mathbb R^3}V(\tau)\Big(\int_0^{\infty}
u^{-2n-3}e^{-if/u}\,du\Big)\,d\tau.
\end{split}
\end{equation*}
We first look at the inner integral:
$$
\int_0^{\infty} u^{-2n-3}e^{-if/u}\,du.
$$
Changing variable $v=\frac{if}{u}$, $u=\frac{if}{v}$, yields
$dv=\frac{-if}{u^2}\,du$ and $du=-\frac{if}{v^2}\,dv$. Hence, we
have
$$
\int_0^{\infty} u^{-2n-3}e^{-if/u}\,du=\frac{1}{i^{2n+2}f^{2n+2}}
\int_0^{\infty}
e^{-v}v^{2n+3-2}\,dv=\frac{\Gamma(2n+2)}{i^{2n+2}f^{2n+2}}.
$$
Let us introduce the following notation
$$
-G(x,z)=\int_0^{\infty} P(x,z,u)\,du
=C\frac{\Gamma(2n+2)}{i^{2n+2}}\int_{\mathbb R^3}
\frac{V(\tau)}{f^{2n+2}(x,z,\tau)}\,d\tau.
$$
The aim of this section is to show that the function $-G(x,z)$ is
the Green function for the sub-Laplacian operator. Firstly, we
need some auxiliary results.

\begin{proposition}
\label{prop:est} Denote
\[
f(x,z,w)=\sum_{l=1}^{n}\frac{1}{4}|x_l|^2|w|_l\coth(
|w|_l)-i\sum_{m=1}^3 w_mz_m=\gamma(x,w)-i\sum_{m=1}^3 w_mz_m,
\] where $w=\tau+i\varepsilon \tilde{z}$.
Then there exist positive constants $c_1$, $c_2$, and
$\varepsilon_0$ such that for all real $\tau_m$, all
$0<\varepsilon<\varepsilon_0$, and all $x\in {\mathbb R}^{4n}$,
$z=(z_1,z_2,z_3)\in {\mathbb R}^3$ we have the estimates
\begin{equation}\label{h15}
\big|\im(\gamma)(x,\tau+i\varepsilon \tilde{z})\big|\le
c_1\varepsilon |x|^2,\end{equation}
\begin{equation}\label{h16}
\re(\gamma)(x,\tau+i\varepsilon \tilde{z})\ge
c_2|x|^2,\end{equation}
\begin{equation}\label{h17}\re(f)(x,z,\tau+i\varepsilon \tilde{z})\ge
c_2(|x|^2+\varepsilon|z|).
\end{equation}
Here $\tilde{z}=\frac{z}{|z|}$ if $z\ne 0$ and $\tilde{z}=0$ if
$z=0$.
\end{proposition}
\begin{proof}
If $\tilde{z}=0$, then $\im (\gamma)(x,\tau)=0$ and since
$|\tau|_l\coth(|\tau|_l)\geq 1$, $l=1,\ldots,n$, we have
$\re(\gamma)(x,\tau)\ge \frac{|x|^2}{4}$.

Suppose that $\tilde{z}\neq 0$. We denote $|w|_l=\Big(\sum_{m=1}^3
a^2_{ml}\big(\tau_m+i\varepsilon
\tilde{z}_m\big)^2\Big)^{1/2}=\alpha_l+i\beta_l$, where
$$\alpha_l=\Big(\big(|\tau|_l^2-\varepsilon^2|\tilde z|_l^2\big)^2
+\big(2\varepsilon\sum_m
a^2_{ml}\tau_m\tilde{z}_m\big)^2\Big)^{1/4}
\cos\frac{\arctan\Big(\frac{2\varepsilon\sum_m
a^2_{ml}\tau_m\tilde{z}_m}{|\tau|_l^2-\varepsilon^2|\tilde
z|_l^2}\Big)}{2}+\pi d,\quad d=0,1,$$ and
$$\beta_l=\Big(\big(|\tau|_l^2-\varepsilon^2|\tilde
z|_l^2\big)^2+\big(2\varepsilon\sum_m
a^2_{ml}\tau_m\tilde{z}_m\big)^2\Big)^{1/4}
\sin\frac{\arctan\Big(\frac{2\varepsilon\sum_m
a^2_{ml}\tau_m\tilde{z}_m}{|\tau|_l^2-\varepsilon^2|\tilde
z|_l^2}\Big)}{2}+\pi d,\quad d=0,1.$$ We consider the case $d=0$,
another one can be treated similarly. Since
$$\coth(\alpha+i\beta) =  \frac{\sinh 2\alpha}{\cosh
2\alpha-\cos 2\beta}- i\frac{\sin 2\beta}{\cosh 2\alpha-\cos
2\beta},$$ we have
\begin{equation*}
\begin{split} \re(\gamma)(x,\tau+i\varepsilon \tilde{z}) +
 i\im(\gamma)(x,\tau+i\varepsilon \tilde{z})=&
\sum_{l=1}^{n}\frac{|x_l|^2}{4}(\alpha_l+i\beta_l)\coth(\alpha_l+i\beta_l)\\
= &\sum_{l=1}^{n}
\frac{|x_l|^2}{4}\Big(\frac{\alpha_l\sinh\alpha_l\cosh\alpha_l
+\beta_l\sin\beta_l\cos\beta_l}{\sinh^2\alpha_l+\sin^2\beta_l}\Big)\\
+ &i
\sum_{l=1}^{n}\frac{|x_l|^2}{4}\Big(\frac{\beta_l\sinh\alpha_l\cosh\alpha_l
-\alpha_l\sin\beta_l\cos\beta_l}{\sinh^2\alpha_l+\sin^2\beta_l}\Big).
\end{split}
\end{equation*}

Denotes by $\psi_l$ the angle between nonzero vectors
$(a_{1l}\tau_1,a_{2l}\tau_2,a_{3l}\tau_3)$ and $(a_{1l}\tilde
z_1,a_{2l}\tilde z_2,a_{3l}\tilde z_3)$. We consider two cases,
when $\cos\psi_l=0$ for all $l=1,\ldots,n$, and
$\cos\psi_l=\vartheta_l\neq 0$ for some index~$l$.

\noindent{\bf Case 1.} If $\cos\psi_l=0$ for all $l=1,\ldots,n$,
then $\sum_{m=1}^{3} a^2_{ml}\tau_m\tilde{z}_m=0$. We have
$$\alpha_l=\big(|\tau|_l^2-\varepsilon^2|\tilde
z|_l^2\big)^{1/2}\quad\text{and}\quad \beta_l=0.$$ It gives
$$\im(\gamma)(x,\tau+i\varepsilon
\tilde{z})=\sum_{l=1}^{n}\frac{|x_l|^2}{4}\Big(\frac{\beta_l\sinh\alpha_l\cosh\alpha_l
-\alpha_l\sin\beta_l\cos\beta_l}{\sinh^2\alpha_l+\sin^2\beta_l}\Big)
=0,$$
$$\re(\gamma)(x,\tau+i\varepsilon \tilde{z}) =\sum_{l=1}^{n}\frac{|x_l|^2}{4}\alpha_l\coth(\alpha_l)\geq
\frac{|x|^2}{4},\qquad\forall\alpha_l\in\mathbb R,$$ because
$\alpha_l\coth(\alpha_l)\geq 1$.

\noindent{\bf Case 2.} If $\cos\psi_l=\vartheta_l\neq 0$ for some
$l=1,\ldots,n$, then $\sum_{m=1}^{3}
a^2_{ml}\tau_m\tilde{z}_m=\vartheta_l|\tau|_l|\tilde{z}|_l$. We
can suppose that $\varepsilon$ satisfies
$0<\varepsilon^2<\min_{l=1,\ldots,n}\big\{\frac{|\tau|_l^2}{2|\tilde{z}|_l^2}\big\}$.
We obtain
\begin{eqnarray}\label{h18}
& \frac{2\varepsilon\vartheta_l|\tilde{z}|_l}{|\tau|_l}
<\frac{2\varepsilon\vartheta_l|\tau|_l|\tilde{z}|_l}{|\tau|_l^2-\varepsilon^2|\tilde{z}|_l^2}
<\frac{4\varepsilon\vartheta_l|\tilde{z}|_l}{|\tau|_l},\nonumber
\\ &
k_1|\tau|_l<\Big(\frac{|\tau|_l^4}{4}+\big(2\varepsilon\vartheta_l|\tau|_l|\tilde{z}|_l\big)^2\Big)^{1/4}
<\Big(\big(|\tau|_l^2-\varepsilon^2|\tilde{z}|_l^2\big)^2+\big(2\varepsilon\sum_m
a^2_{ml}\tau_m\tilde{z}_m\big)^2\Big)^{1/4}
\\ & <\Big(|\tau|_l^4+\big(2\varepsilon\vartheta_l|\tau|_l|\tilde{z}|_l\big)^2\Big)^{1/4}
<k_2(\vartheta)|\tau|_l.\nonumber\end{eqnarray} Now, we put one
more restriction to $\varepsilon$ assuming that
$\varepsilon<\min_{l=1,\ldots,n}\big\{\frac{|\tau|_l}{4\vartheta_l|\tilde{z}|_l}\big\}$.
Then
$$\frac{2\varepsilon\vartheta_l|\tilde{z}|_l}{\pi|\tau|_l}
\leq\frac{1}{2}\arctan\frac{2\varepsilon\vartheta_l|\tilde{z}|_l}{|\tau|_l}
\leq\frac{1}{2}\arctan\frac{2\varepsilon\vartheta_l|\tau|_l|\tilde{z}|_l}{|\tau|_l^2-\varepsilon_l^2|\tilde{z}|_l^2}
\leq\frac{1}{2}\arctan\frac{4\varepsilon\vartheta_l|\tilde{z}|_l}{|\tau|_l}
\leq\frac{2\varepsilon\vartheta_l|\tilde{z}|_l}{|\tau|_l},$$ and
we get
\begin{eqnarray}\label{h19}&
\frac{\sqrt 2}{2}
<\cos\Big(\frac{2\varepsilon\vartheta_l|\tilde{z}|_l}{|\tau|_l}\Big)
<\cos\frac{\arctan\big(\frac{2\varepsilon\vartheta_l|\tau|_l|\tilde{z}|_l}{|\tau|_l^2-\varepsilon^2|\tilde{z}|_l^2}\big)}{2}
<\cos\Big(\frac{2\varepsilon\vartheta_l|\tilde{z}|_l}{\pi|\tau|_l}\Big)<1,
\\ & \frac{4\varepsilon\vartheta_l|\tilde{z}|_l}{\pi^2|\tau|_l}<
\sin\Big(\frac{2\varepsilon\vartheta_l|\tilde{z}|_l}{\pi|\tau|_l}\Big)
<\sin\frac{\arctan\big(\frac{2\varepsilon\vartheta_l|\tau|_l|\tilde{z}|_l}{|\tau|_l^2-\varepsilon^2|\tilde{z}|_l^2}\big)}{2}
<\sin\Big(\frac{2\varepsilon\vartheta_l|\tilde{z}|_l}{|\tau|_l}\Big)
<\frac{2\varepsilon\vartheta_l|\tilde{z}|_l}{|\tau|_l}.\nonumber
\end{eqnarray} We observe that $|\tilde{z}|_l^2=\sum_{m=1}^{3}a^2_{ml}\frac{z_m^2}{|z|^2}\leq \sum_{m=1}^{3}a^2_{ml}
\leq \overline a$, where $\overline
a=\max\limits_{m,l}\{a^2_{ml}\}$. From the other hand, if we
denote $\underline a=\min\limits_{m,l}\{a^2_{ml}\}$, then
$\underline a\leq
\sum_{m=1}^{3}a^2_{ml}\frac{z_m^2}{|z|^2}=|\tilde{z}|_l^2$.
From~\eqref{h18} and~\eqref{h19} we estimate the value of
$\alpha_l$ and $\beta_l$ as follows
\begin{equation}\label{h13}
k_1|\tau|_l<\alpha_l<k_2|\tau|_l,\end{equation}\begin{equation}\label{h14}
k_3(\underline a)\varepsilon<\beta_l<k_4(\overline a)\varepsilon.
\end{equation}
If $|\tau|_l< 1$ we use the Taylor decomposition and obtain
$$\Big|\frac{\beta_l\sinh\alpha_l\cosh\alpha_l-\alpha_l\sin\beta_l\cos\beta_l}{\sinh^2\alpha_l+\sin^2\beta_l}\Big|
=\Big|\frac{-\frac{2}{3}\alpha_l\beta_l(\alpha_l^2-\beta_l^2)+O(\alpha_l^4-\beta_l^4)}{\alpha_l^2+\beta_l^2
-O(\alpha_l^4+\beta_l^4)}\Big|\leq k_5\varepsilon.$$ If
$|\tau|_l\geq 1$ we argue as follows
$$\Big|\frac{\beta_l\sinh\alpha_l\cosh\alpha_l-\alpha_l\sin\beta_l\cos\beta_l}{\sinh^2\alpha_l+\sin^2\beta_l}\Big|
\leq|\beta_l|\big(|\coth(\alpha_l)|+\big|\frac{\alpha_l}{\sinh^2(\alpha_l)}\big|\big)\leq
k_6\varepsilon,$$ because $\alpha_l$ is bounded from below, the
functions $|\coth(\alpha_l)|$ and
$\big|\frac{\alpha_l}{\sinh^2\alpha_l}\big|$ are bounded from above.
The last two estimates imply
$$|\im(\gamma)(x,\tau+i\varepsilon \tilde{z})|\leq \sum_{n}^{l=1}\frac{|x|_l^2}{4}k_7\varepsilon\leq
c_1\varepsilon|x|^2.$$

To obtain~\eqref{h16} we change the arguments. Let us focus on the
value of the derivatives $\frac{\partial \gamma(x,w)}{\partial w_m}$
at $w_m=i\zeta_m$, $\zeta_m\in\mathbb R$ for $m=1,2,3$. The equality
$$\frac{\partial \gamma(x,w)}{\partial
w_m}\Big|_{w=i\zeta}=-\sum_{l=1}^{n}\frac{|x|_l^2}{4}\frac{a^2_{ml}w_m}{|w|_l}\Big(\frac{|w|_l}{\sinh^2(|w|_l)}-\coth(|w|_l)\Big)
\Big|_{w=i\zeta}
=i\sum_{l=1}^{n}\frac{|x|_l^2}{4}\frac{a^2_{ml}\zeta_m}{|\zeta|_l}\mu(|\zeta|_l)$$
implies $\frac{\partial \re(\gamma(x,w))}{\partial
w_m}\Big|_{w=i\zeta}=0$ and we conclude that $w=i\zeta$ is a
critical point for $\re(\gamma(x,w))$. Let us look at the Hessian
at $w=i\zeta$. We have
$$\frac{\partial^2 \gamma}{\partial w_m^2}\Big|_{w=i\zeta}=\sum_{l=1}^{n}\frac{|x|_l^2a^2_{ml}}{4}
\Big[\frac{\mu(|\zeta|_l)}{|\zeta|_l}\big(1-\frac{a^2_{ml}\zeta_m^2}{|\zeta|_l^2}\big)
+\frac{2a^2_{ml}\zeta_m^2}{|\zeta|_l^2\sin^2(|\zeta|_l)}\big(1-|\zeta|_l\cot(|\zeta|_l)\big)\Big].$$
Since $1-\frac{a^2_{ml}\zeta_m^2}{|\zeta|_l^2}\geq 0$ and
$1-|\zeta|_l\cot(|\zeta|_l)\geq 0$ we see that $\frac{\partial^2
\gamma}{\partial w_m^2}\Big|_{w=i\zeta}>0$ for $0\neq
|\zeta|_l<\frac{\pi}{2}$, $m=1,2,3$. The mixed second derivatives
are
\begin{eqnarray*}\frac{\partial^2 \gamma}{\partial w_m\partial
w_k}\Big|_{w=i\zeta} & = &
\sum_{l=1}^{n}\frac{|x_l|^2}{4}\frac{a^2_{ml}a^2_{kl}\zeta_m\zeta_k}{|\zeta|_l^2}
\Big[\frac{\cot(|\zeta|_l)}{|\zeta|_l}-\frac{2|\zeta|_l\cot(|\zeta|_l)}{\sin^2(|\zeta|_l)}
+\frac{1}{\sin^2(|\zeta|_l)}\Big] \\ & = &
\sum_{l=1}^{n}\frac{|x_l|^2}{4}\frac{a^2_{ml}a^2_{kl}\zeta_m\zeta_k}{|\zeta|_l^2}g(|\zeta|_l),\end{eqnarray*}
where
$g(|\zeta|_l)=\frac{\cot(|\zeta|_l)}{|\zeta|_l}-\frac{2|\zeta|_l\cot(|\zeta|_l)}{\sin^2(|\zeta|_l)}
+\frac{1}{\sin^2(|\zeta|_l)}$. We observe that since all second
derivatives of $\gamma(x,w)$ are real at the critical point, the
Hessian for $\gamma(x,w)$ coincides with the Hessian $H$ for
$\re(\gamma(x,w))$ at $w=i\zeta$. We write $\frac{\partial^2
\gamma}{\partial w_m^2}$ as
$$\frac{\partial^2 \gamma}{\partial
w_m^2}\Big|_{w=i\zeta}=\sum_{l=1}^{n}\frac{|x|_l^2}{4}\Big[\frac{a^2_{ml}\mu(|\zeta|_l)}{|\zeta|_l}+
\frac{a^4_{ml}\zeta_m^2}{|\zeta|_l^2}g(|\zeta|_l)\Big].$$ Then the
Hessian can be written in the form $H=\sum_{l=1}^{n}H_l$, where
\begin{equation*}
H_l=\frac{|x|_l^2}{4}\left[\array{ccc}
\frac{a^2_{1l}\mu(|\zeta|_l)}{|\zeta|_l}+
\frac{a^4_{1l}\zeta_1^2}{|\zeta|_l^2}g(|\zeta|_l) &
\frac{a^2_{1l}a^2_{2l}\zeta_1\zeta_2}{|\zeta|_l^2}g(|\zeta|_l) &
\frac{a^2_{1l}a^2_{3l}\zeta_1\zeta_3}{|\zeta|_l^2}g(|\zeta|_l)
\\
\frac{a^2_{1l}a^2_{2l}\zeta_1\zeta_2}{|\zeta|_l^2}g(|\zeta|_l) &
\frac{a^2_{2l}\mu(|\zeta|_l)}{|\zeta|_l}+
\frac{a^4_{2l}\zeta_2^2}{|\zeta|_l^2}g(|\zeta|_l) &
\frac{a^2_{2l}a^2_{3l}\zeta_2\zeta_3}{|\zeta|_l^2}g(|\zeta|_l)
\\
\frac{a^2_{1l}a^2_{3l}\zeta_1\zeta_3}{|\zeta|_l^2}g(|\zeta|_l) &
\frac{a^2_{2l}a^2_{3l}\zeta_2\zeta_3}{|\zeta|_l^2}g(|\zeta|_l) &
\frac{a^2_{3l}\mu(|\zeta|_l)}{|\zeta|_l}+
\frac{a^4_{3l}\zeta_3^2}{|\zeta|_l^2}g(|\zeta|_l)
\endarray\right].\end{equation*} To show that $H$ is positive definite we need to show that each $H_l$
is positive definite. It was shown that
$$\frac{a^2_{1l}\mu(|\zeta|_l)}{|\zeta|_l}+
\frac{a^4_{1l}\zeta_1^2}{|\zeta|_l^2}g(|\zeta|_l)>0,\quad\text{for}\quad
0\neq |\zeta|_l<\frac{\pi}{2}.$$ Then we have
\begin{equation*}
\begin{split}
&\Big(\frac{a^2_{1l}\mu(|\zeta|_l)}{|\zeta|_l}+
\frac{a^4_{1l}\zeta_1^2}{|\zeta|_l^2}g(|\zeta|_l)\Big)\Big(\frac{a^2_{2l}\mu(|\zeta|_l)}{|\zeta|_l}+
\frac{a^4_{2l}\zeta_2^2}{|\zeta|_l^2}g(|\zeta|_l) \Big)-\Big(\frac{a^2_{1l}a^2_{2l}\zeta_1\zeta_2}{|\zeta|_l^2}g(|\zeta|_l)\Big)^2\\
& =
\frac{a^2_{1l}a^2_{2l}\mu^2(|\zeta|_l)}{|\zeta|_l^2}\Big(1-\frac{a^2_{1l}\zeta_1^2+a^2_{2l}\zeta_2^2}{|\zeta|_l^2}\Big)+
\frac{2a^2_{1l}a^2_{2l}\mu(|\zeta|_l)(a^2_{1l}\zeta_1^2+a^2_{2l}\zeta_2^2)}{|\zeta|_l^4\sin^2(|\zeta|_l)}\big(1-|\zeta|_l\cot(|\zeta|_l)\big)
>0\end{split}
\end{equation*} for
$0\neq |\zeta|_l<\frac{\pi}{2}.$ Finally, we calculate $\det H_l$:
\begin{equation*}
\frac{|x|_l^6}{4^3}\frac{a^2_{1l}a^2_{2l}a^2_{3l}\mu^2(|\zeta|_l)}{|\zeta|_l^2}
\Big(\frac{\mu(|\zeta|_l)}{|\zeta|_l}+g(|\zeta|_l)\Big)=
\frac{a^2_{1l}a^2_{2l}a^2_{3l}|x|_l^6}{4^3}\frac{2\mu^2(|\zeta|_l)}{|\zeta|_l^2\sin^2(|\zeta|_l)}
\big(1-|\zeta|_l\cot(|\zeta|_l)\big)>0
\end{equation*} for $0\neq|\zeta|_l<\frac{\pi}{2}$. We conclude that the Hessian is positive definite and
$\re(\gamma(z,w))$ has a local minimum at $w=i\zeta$. Thus
$$\re(\gamma(x,w))\geq \re(\gamma(x,w))\vert_{w=i\zeta}=
\sum_{l=1}^{n}\frac{|x|_l^2}{4}|\zeta|_l\cot(|\zeta|_l)\geq
c_2|x|^2\quad\text{if} \quad |\zeta|_l<\pi/4.$$ Put
$\zeta=\varepsilon\tilde z$, then
$|\zeta|_l\leq\varepsilon\overline a$ and~\eqref{h16} holds with
$\varepsilon_0=\frac{\pi}{4\overline a}$.

Estimate~\eqref{h17} is a consequence of estimates~\eqref{h15}
and~\eqref{h16} since
\begin{eqnarray*}
f(x,z,\tau+i\varepsilon\tilde z) &
=\gamma(x,z,\tau+i\varepsilon\tilde
z)-i\sum\limits_{m=1}^{3}(\tau_m+i\varepsilon\frac{z_m}{|z|})z_m\\
& =\gamma(x,z,\tau+i\varepsilon\tilde
z)+\varepsilon|z|-i\sum\limits_{m=1}^{3}\tau_mz_m.
\end{eqnarray*}
\end{proof}

Theorem~\ref{prop:est} in a non-diagonal situation was proved
in~\cite{BGGR3}.
\medskip
\begin{lemma}
\label{lem:fund1} If $x$ is a non-zero vector in $\R^{4n}$, the
integral
\begin{equation}
\label{eq:fund3} \tilde G(x,z)=\int_{{\mathbb
R}^3}\frac{V(\tau)}{f^{2n+2}(x,z,\tau)}d\tau
\end{equation}
is absolutely convergent and one has for $x\ne 0$
\[
\Delta_0\tilde G(x,z)=0.
\]
\end{lemma}
\begin{proof}
Since the function $V(\tau)$ does not depend on $x$ and $z$, we
have $\Delta_0 V=0$, $X_{kl}V=0$, $k=1,2,3,4$, $l=1,\ldots,n$, and
the equation~\eqref{an12} reduced to the following one
\begin{equation}\label{an13}\Delta_0(Vf^{-p})=-p\Big(f^{-p-1}V\Delta_0f+(-p-1)Vf^{-p-2}H(Xf)\Big).\end{equation}
Here $H(Xf)=|Xf|^2$ by Proposition~\ref{h5}. Moreover, taking into
account that the complex action function $f(x,z,\tau)$ satisfies
the Hamilton-Jacobi equation~\eqref{eq:HJ} and $p=2n+2$, we get
\begin{equation}\label{an14}
\Delta_0(Vf^{-2n-2})=(-2n-2)\Big(f^{-2n-3}V\Delta_0f+(-2n-3)Vf^{-2n-3}
+(2n+3)Vf^{-2n-4}\sum_{m=1}^{3}\tau_m\frac{\partial f}{\partial
\tau_m}\Big).\end{equation} Substituting the last term in the
right hand side of~\eqref{an14} from the formula
$$\sum_{m=1}^{3}\frac{\partial}{\partial
\tau_m}\big(\tau_mVf^{-2n-3}\big)=3Vf^{-2n-3}+f^{-2n-3}\sum_{m=1}^{3}\tau_m\frac{\partial
V}{\partial
\tau_m}-(2n+3)Vf^{-2n-4}\sum_{m=1}^{3}\tau_m\frac{\partial
f}{\partial \tau_m},$$ we deduce
$$\Delta_0(Vf^{-2n-2})=(-2n-2)\Big(f^{-2n-3}\big(V(\Delta_0f-2n)+\sum_{m=1}^{3}\tau_m\frac{\partial
V}{\partial \tau_m}\big)-\sum_{m=1}^{3}\frac{\partial}{\partial
\tau_m}\big(\tau_mVf^{-2n-3}\big)\Big).$$ Since the volume element
$V(\tau)$ is a solution of the transport equation~\eqref{h10},
finally, we obtain $$\Delta_0\tilde G(x,z)=\int_{\mathbb
R^3}\Delta_0(Vf^{-2n-2})\,d\tau=(2n+2)\int_{\mathbb
R^3}\sum_{m=1}^{3}\frac{\partial}{\partial
\tau_m}\big(\tau_mVf^{-2n-3}\big)\,d\tau.$$ We observe that
\begin{equation}\label{an15}V(\tau)=\prod_{l=1}^{n}\frac{|\tau|_l^{2}}{\sinh^{2}|\tau|_l}\to
0\quad \text{as one of the}\quad |\tau_m|=R\to
\infty,\end{equation} and
\begin{equation}\label{an16}
|f|\geq\sum\limits_{l=1}^{n}\frac{|x_l|^2}{2}|\tau|_l\coth(|\tau|_l)\geq\frac{|x|^2}{2}
\end{equation} because of $|\tau|_l\coth(|\tau|_l)\geq 1$ for $l=1,\ldots,n$.
The estimates~\eqref{an15} and~\eqref{an16} show
\begin{equation}\label{an34}\lim_{R\to\infty}\int_{|\tau_m|\leq
R}\frac{\partial}{\partial
\tau_m}\big(\tau_mVf^{-2n-3}\big)=0.\end{equation} The last
equality implies
$$\Delta_0(Vf^{-2n-2})=(2n+2)\int_{\mathbb
R^3}\sum_{m=1}^{3}\frac{\partial}{\partial
\tau_m}\big(\tau_mVf^{-2n-3}\big)\,d\tau=0,$$ that terminates the
proof of Lemma~\ref{lem:fund1}.
\end{proof}

The argument of Lemma \ref{lem:fund1} is not valid for $x=0$. In
fact, the integral (\ref{eq:fund3}) is divergent because the
denominator of the integrand contains $-\sum_{m=1}^3 \tau_mz_m$
which is zero along a hyperplane of ${\mathbb R}^3$. We will treat
the case by changing contour by adding a small imaginary part to
the $\tau_m$'s. We shall prove that for $x\ne 0$ we can change the
contour in (\ref{eq:fund3}) and when $x$ will be zero the integral
(\ref{eq:fund3}) will still be convergent on the new contour. In
order to achieve this goal, we need to use Proposition
\ref{prop:est}.
\begin{proposition} For $x\ne 0$, the integral $\tilde G(x,z)$ defined
in (\ref{eq:fund3}) is given by
\begin{equation}
\label{eq:fund4} \tilde G(x,z)=\int_{{\mathbb
R}^3}\frac{V(\tau+i\varepsilon
\tilde{z})}{f^{2n+2}(x,z,\tau+i\varepsilon \tilde {z})}d\tau
\end{equation}
for $0<\varepsilon<\varepsilon_0$ sufficiently small. The integral
(\ref{eq:fund4}) makes sense even for $x=0$ and $z\ne 0$, so the
function $\tilde G(x,z)$ is well-defined (in fact, real analytic)
except at the origin in ${\mathbb R}^{4n}\times {\mathbb R}^3$ and
satisfies
\[
\Delta_0\tilde G(x,z)=0\qquad{\hbox{for}}\quad (x,z)\ne (0,0).
\]
\end{proposition}
\begin{proof} We may prove this theorem by imitating the idea in \cite{BGGR3}.
Set
\[
\Omega_{K,\varepsilon}=\{\xi=\tau+i\eta\tilde
z:\,\tau\in\R^3,\,\,|\tau|<K,\,\,0<\eta<\varepsilon\}\,\subset\,{\mathbb
C}^3.
\]
Assume that $|x|\ne 0$. By Theorem \ref{th:heat1} the differential
form
\[ \omega=\big(V(\zeta)/f^{2n+2}(x,z,\zeta)\big)d\xi_1\wedge
d\xi_2\wedge d\xi_3 \] is a homomorphic form of type $(3, 0)$ in
$\Omega_{K,\varepsilon}$. It is easy to see that its differential is
zero. Hence, by Stokes's Theorem
\[
\int_{\partial \Omega_{K,\varepsilon}}\omega=0.
\]
The boundary can be written as $\partial \Omega_{K,\varepsilon}
=\partial \Omega_1\cup\partial \Omega_2\cup\partial \Omega_3$.
\begin{itemize} \item[(1)] {The set $\partial
\Omega_1=\{\tau\in\R^3:\,|\tau|<K,\,\eta=0\}$ which is such that
\[
\lim_{K\rightarrow\infty}\int_{|\tau|<K}\frac{V(\tau)}{f^{2n+2}(x,z,\tau)}d\tau=\tilde
G(x,z)
\]
since the integral (\ref{eq:fund3}) converges absolutely.}
\item[(2)] {The set $\partial \Omega_2=\{\tau+i\varepsilon\tilde
z:\,|\tau|<K,\,\eta=\varepsilon\}$. The integral (\ref{eq:fund4})
converges absolutely by Proposition \ref{prop:est} in this case.}
\item[(3)] {The set $\partial
\Omega_3=\{\xi=\tau+i\varepsilon\tilde
z:\,|\tau|=K,\,\,0<\eta<\varepsilon\}$. Again, by Proposition
\ref{prop:est}, one has
\[
\lim_{K\rightarrow\infty}\int_{\partial
\Omega_3}\frac{V(\tau+i\eta \tilde{z})}{f^{2n+2}(x,z,\tau+i\eta
\tilde {z})}d\xi=0.
\]}\end{itemize}

By the discussion above, one may conclude that for $x\ne 0$,
$\tilde G(x,z)$ is given by the integral (\ref{eq:fund4}) on a
shifted contour. Moreover, this integral is absolutely convergent
even when $x=0$ and $z\ne 0$ by Proposition \ref{prop:est}. We
complete the proof of this theorem.
\end{proof}
\begin{theorem}
\label{th:fund5} The kernel $G(x,z)$ of the Green's function for
the sub-Laplacian $\Delta_0$ is given by the formula
\[
G(x,z)=-\frac{2^{2n}(2\pi)^{2n+3}}{(2n+1)!}\int_{{\mathbb
R}^3}\frac{V(\tau+i\varepsilon
\tilde{z})}{f^{2n+2}(\tau+i\varepsilon \tilde{z})}d\tau.
\]
\end{theorem}
\begin{proof} For any $K>0$, denote
\[
\tilde G_K(x,z)=\frac{1}{\Gamma(2n+2)}\int_{{\mathbb
R}^3}V(\tau+i\varepsilon \tilde{z})d\tau\int_0^K
t^{2n+2-1}e^{-tf}dt.
\]
The function $\tilde G_K(x,z)$ is smooth everywhere and for
$(x,z)\ne (0,0)$, one has
\[
\lim_{K\rightarrow \infty}\tilde G_K(x,z)=\tilde G(x,z).
\]
Using Proposition \ref{prop:est}, we know that
\[
\big|\tilde G_K(x,z)\big|\le C\int_{{\mathbb
R}^3}d\tau\int_0^\infty
\prod_{l=1}^{n}\frac{|\tau|_l^2}{\sinh^2|\tau|_l}\cdot
e^{-c_1(|x|^2+|z|)t}t^{2n+1}dt \le C(|x|^2+|z|)^{-2n-2}.
\]
But the function $(|x|^2+|z|)^{-2n-2}\in L^1_{loc}({\mathbb
R}^{4n+3})$ since the homogeneous degree is $4n+6$ in this case.
It follows that
\[
\lim_{K\rightarrow \infty}\tilde G_K(x,z)=\tilde G(x,z)\qquad
{\hbox {in}}\quad L^1_{loc}({\mathbb R}^{4n+3})
\]
by the Dominated Convergence Theorem. We first calculate
$\Delta_0\tilde G_R(x,z)$ for $x\ne 0$. We need to compute
$\Delta_0(Ve^{-tf})$ which is
\[
\Delta_0(Ve^{-tf})=-tV\cdot\Delta_0 (f)e^{-tf}+t^2H(x,z,\nabla
f)Ve^{-tf}.
\]
But the Hamilton-Jacobi equation~\eqref{eq:HJ} yields
\begin{equation*}
\begin{split}
\Delta_0(Ve^{-tf})=& -tV\cdot \Delta_0(f)e^{-tf}+t^2\Big(f-
\sum_{m=1}^3 \tau_m\frac{\partial f}{\partial
\tau_m}\Big)Ve^{-tf}\\
= & -te^{-tf}\Big(V\Delta_0(f)+\sum_{m=1}^3
\frac{\partial}{\partial \tau_m}(\tau_mV)\Big)+t^2e^{-tf}Vf\Big)\\
&+t\sum_{m=1}^3\frac{\partial}{\partial
\tau_m}\big(\tau_me^{-tf}V\big).
\end{split}
\end{equation*}
We also know that
\[
\int_0^K t^{2n+3}e^{-tf}fdt=-K^{2n+3}e^{-Kf}+(2n+3)\int_0^K
t^{2n+2}e^{-tf}dt.
\]
This implies that
\begin{equation}\label{an35}
\begin{split}
\int_{\R^3}d\tau\int_0^K
t^{2n+1}e^{-tf}t^2fVdt=&-K^{2n+3}\int_{\R^3}Ve^{-Kf}d\tau\\
&+(2n+3)\int_0^Kt^{2n+2}dt\int_{\R^3}e^{-tf}Vd\tau.
\end{split}
\end{equation}
Hence,
\begin{equation*}
\begin{split}
\Delta_0\tilde G_K(x,z)=&
\frac{-1}{\Gamma(2n+2)}\int_0^Kt^{2n+2}dt\int_{\R^3}e^{-tf}
\Big[\sum_{m=1}^3 \frac{\partial(\tau_mV)}{\partial
\tau_m}+(\Delta_0 f-2n-3)V\Big]d\tau\\
&-\frac{K^{2n+3}}{\Gamma(2n+2)}\int_{\R^3}Ve^{-Kf}d\tau,
\end{split}
\end{equation*} where we used~\eqref{an34} and~\eqref{an35}.
The first integral vanishes since $V$ satisfies the generalized
transport equation~\eqref{h10}. Therefore, for $x\ne 0$,
\[
\Delta_0\tilde
G_K(x,z)=-\frac{K^{2n+3}}{\Gamma(2n+2)}\int_{\R^3}V(\tau)e^{-Kf(x,z,\tau)}d\tau.
\]
However, we may also change contour in the above integral and
obtain
\begin{equation}
\label{eq:fund6} \Delta_0\tilde
G_K(x,z)=-\frac{K^{2n+3}}{\Gamma(2n+2)}\int_{\R^3}V(\tau+i\varepsilon
\tilde{z})e^{-Kf(x,z,\tau+i\varepsilon
\tilde{z})}d\tau,\quad{\hbox{for}}\quad x\ne 0.
\end{equation}
Since $\tilde G_K(x,z)$ is smooth everywhere, the integral
(\ref{eq:fund6}) provides the value of $\Delta_0\tilde G_K(x,z)$ at
every point where the integral is convergent. It follows that
$\Delta_0\tilde G_K(x,z)$ is equal to
\[-\frac{K^{2n+3}}{(2n+1)!}\int_{\R^3}V(\tau+i\varepsilon
\tilde{z})e^{-Kf(x,z,\tau+i\varepsilon \tilde{z})}d\tau\]
everywhere and to
$-\frac{K^{2n+3}}{(2n+1)!}\int_{\R^3}V(\tau)e^{-Kf(x,z,\tau)}d\tau$
almost everywhere. Furthermore, for $(x,z)$ in a compact set $U$
with $U$ disjoint from the origin,
\[
\re(f)(x,z,\tau+i\varepsilon \tilde{z})\ge \kappa>0.
\]
Hence $\Delta_0\tilde G_K(x,z)\,\rightarrow\,0$ uniformly as
$K\rightarrow\infty$ on compact subsets of $\R^{4n}\times \R^3$
which is disjoint from the origin. Now we need to compute the
$L^1$-norm of $\Delta_0\tilde G_K(x,z)$. Since the integral
$\int_{\R^3}V(\tau)e^{-tf(x,z,\tau)}d\tau$ coincides with the
integral of the right-hand side of (\ref{eq:fund6}) almost
everywhere, we can just compute the following
\[
{\mathcal I}:=-\frac{K^{2n+3}}{(2n+1)!}\Big|
\int_{\R^{4n}}\int_{\R^3}\int_{\R^3}V(\tau)e^{-K[\gamma(x,\tau)-i\sum_{m=1}^3\tau_mz_m]}d\tau
dzdx\Big|.
\]
Since the integral converges absolutely, we may interchange the
order of the integration by Fubini's Theorem. Let us integrate the
$z$-variable first.
\begin{equation*}
\begin{split}
&\int_{\R^3}e^{-K[\gamma(x,\tau)-i\sum_{m=1}^3\tau_mz_m]}dz
=\int_{\R^3}
e^{-\frac{K}{4}\sum_{l=1}^{n}|x|_l^2|\tau|_l\coth(|\tau|_l)+iK\sum_{m=1}^3\tau_mz_m]}dz \\
=&e^{-K[1/4\sum_{l=1}^{n}|x|_l^2|\tau|_l\coth(|\tau|_l)}\int_{\R^3}
e^{iK\sum_{m=1}^3\tau_mz_m}dz\\
=&-e^{-K[1/4\sum_{l=1}^{n}|x|_l^2|\tau|_l\coth(|\tau|_l)}K^{-3}
(2\pi)^3{\mathcal F}^{-1}(1) =K^{-3} (2\pi)^3\delta(\tau).
\end{split}
\end{equation*}
Here ${\mathcal F}(g)(\xi)=\int_{\R^3}g(x)e^{-2\pi ix\cdot \xi}dx$
is the Fourier transform of the function $g(x)$. We know that
\[
\lim_{\tau\rightarrow0}\gamma(x,\tau)=\lim_{\tau\rightarrow0}\frac{1}{4}
\sum_{l=1}^{n}|x|_l^2|\tau|_l\coth(|\tau|_l)=\frac{1}{4}|x|^2
\]
and
\[
\lim_{\tau\rightarrow 0}V(\tau)= \lim_{\tau\rightarrow
0}\prod_{l=1}^{n}\frac{|\tau|_l^2}{\sinh^2|\tau|_l}=1.
\]
It follows that
\[
K^{-3}
(2\pi)^3\int_{\R^3}V(\tau)e^{-K\gamma(x,\tau)}\delta(\tau)d\tau
=\frac{(2\pi)^3}{K^3} e^{-\frac{K}{4}|x|^2}
\]
in the sense of distribution. Finally, one has
\[
\frac{1}{K^3}
(2\pi)^3\int_{\R^{4n}}e^{-\frac{K}{4}|x|^2}dx=\frac{(2\pi)^3}{K^3}
2^{4n}\frac{\pi^{2n}}{K^{2n}}.
\]
This gives us
\[
{\mathcal I}:=-\frac{K^{2n+3}}{(2n+1)!}\frac{(2\pi)^3}{K^3}
2^{4n}\frac{\pi^{2n}}{K^{2n}}=-2^{2n}\frac{(2\pi)^{2n+3}}{(2n+1)!}.
\]
This proves that $\Delta_0\tilde G_K(x,z)\rightarrow 0$ uniformly
on compact sets on $\R^{4n}\times \R^3$ disjoint from the origin
with a constant integral over $\R^{4n}\times\R^3$. This means that
when $K\rightarrow\infty$,
\[
\Delta_0\tilde
G_K(x,z)\,\rightarrow\,-2^{2n}\frac{(2\pi)^{2n+3}}{(2n+1)!}\delta_{(0,0)}.
\]
On the other hand, $\tilde G_K(x,z)\rightarrow\tilde G(x,z)$ in
$L^1_{loc}(\mathbb R^{4n+3})$ as $K\to\infty$. Hence,
\[
\Delta_0\tilde G_K(x,z)\,\rightarrow\,\Delta_0 \tilde G(x,z)
\]
in the sense of distribution. Therefore,
\[
\Delta_0\tilde
G(x,z)=-\frac{2^{2n}(2\pi)^{2n+3}}{(2n+1)!}\delta_{(0,0)}.
\]
The proof of the theorem is therefore complete.
\end{proof} The symmetry of homogeneous $\mathbb H$-type groups
allows us to deduce another form of Green's function related to the
homogeneous norm (see, for instance~\cite{CDG}).
\medskip
\section{Estimates of the fundamental solution}

In this section, we discuss sharp estimates for the integral
operator induced by the fundamental solution $G(x,z)$:
\[
{\bf G}(g)(x,z)=G\ast g(x,z)=\int_{Q^n}G(y,w)g((y,w)^{-1}\cdot
(x,z))dydw \] in $L_k^p(Q^n)$ Sobolev spaces, Hardy-Sobolev spaces
$H_k^p(Q^n)$ for $k\in {\bf Z}_+$, $0<p<\infty$ and Lipschitz
spaces $\Lambda_\beta$, $\Gamma_\beta$ with $\beta>0$. We consider
here the sub-Laplacian $\Delta_0$. It is easy to see from the
group law that the operator is homogeneous of degree $-2$ under
that non-isotropic dilation:
\[ \delta_\lambda:\, (x,z)\,\rightarrow\, (\lambda x,\lambda^2 z).
\] Hence, in general, the homogeneous degree of $Q^n$ is $4n+6$.
As in Folland-Stein~\cite{FS} and Koranyi~\cite{Kor1}, we may
define a homogeneous norm
\[
|(x,z)|=\Big(\big(\sum_{j=1}^n\sum_{k=1}^4
x_{jk}^2\big)^2+\sum_{j=1}^3z_j^2\Big)^{\frac{1}{4}}.
\]
Then we may define a pseudo metric by
$\rho(\bx,\by)=|\bx\cdot\by^{-1}|$ where $\bx=(x,z)$ and
$\by=(y,w)$. Then it is easy to see that $\rho$ is equivalent to
the Carnot-Carath\'eodory metric. Using this metric, one may
obtain estimates of $G(x,z)$ in various function spaces. From the
discussion in Section 7, we know that the fundamental solution
$G(x,z)$ with singularity at the origin is a homogeneous kernel of
degree $-4n-4$. In fact, the fundamental solution $G(x,z)$
satisfies the following size estimates:
\[
\big|G(x,z)\big|\le
\frac{C_1|(x,z)|^2}{\big|B((x,z),|(x,z)|)\big|},\qquad
\Big|\frac{\partial^{|\alpha|}G(x,z)}{\partial
x_{11}^{\alpha_{11}}\cdots
\partial x_{4n}^{\alpha_{4n}}}\Big|\le\frac
{C_2|(x,z)|^{2-|\alpha|}}{\big|B((x,z),|(x,z)|)\big|},
\]
and \[
\Big|\frac{\partial^{\beta_1+\beta_2+\beta_3}G(x,z)}{\partial
z_1^{\beta_1}
\partial z_2^{\beta_2}\partial z_3^{\beta_3}}\Big|\le
\frac
{C_3|(x,z)|^{2-2(\beta_1+\beta_2+\beta_3)}}{\big|B((x,z),|(x,z)|)\big|}.
\]
Therefore, it is a locally integrable function. From classical
results, it is easy to see that the operator ${\bf G}$ originally
defined on the Schwartz space ${\mathcal S}(Q^n)$ can be extended to
a bounded operator from $L^p_{loc}(Q^n)$ into $L^p_{loc}(Q^n)$ for
$1\le p\le \infty$ (see \cite{BGGR3} and \cite{FS} ). Moreover,
${\bf G}$ is a smoothing operator. Hence, we may allow to
differentiate the kernel $G(x,z)$. Moreover, for any $l=1,\ldots, n$
and $k=1,\ldots,4$, $X_{kl}{\bf G}$ originally defined on ${\mathcal
S}(Q^n)$ can be extended to a bounded operator from $L^p_{loc}(Q^n)$
into $L^p_{loc}(Q^n)$ for $1\le p\le \infty$. The problem now
reduces to looking at the second derivatives (in the horizontal
directions) of $G(x,z)$. Before we go further, let us recall some
basic definitions and properties of several functions spaces.
\smallskip

\noindent $\bullet$ {\bf Lipschitz spaces} As in \cite{FS}, we
define the space $\Gamma_\beta(Q^n)$ as the set of all bounded
functions $g$ with compact support on $Q^n$ such that
\[
\sup_{(x,z),(y,w)\in Q^n}\Big|g((x,z)\cdot
(y,w)^{-1})-g(x,z)\Big|\le C\cdot |(y,w)|^{\beta},\qquad
0<\beta<1.
\]
When $k<\beta<k+1$ with $k=1,2,...$, we may define $g\in
\Gamma_\beta(Q^n)$ as the set of all $C^k$ functions $g$ with
compact support on $Q^n$ such that ${\mathcal P}(X,X^\prime)g\in
\Gamma_\beta(Q^n)$, where ${\mathcal P}(X,X^\prime)$ is a monomial
of degree $k$ in vector fields $X,\,X^\prime\in V_1$. For integral
value $k$, one may define the space $\Gamma_\beta(Q^n)$ by
interpolation. Furthermore, since $|(x,z)|\le A\|(x,z)\|^{1/2}$
for $|(x,z)|$ small, one may conclude that
$\Gamma_{\beta}(Q^n)\subset \Lambda_{\beta/2}(Q^n)$. Here
$\|(x,z)\|$ is the Euclidean distance between the point $(x,z)$
and the origin and $\Lambda_\alpha(Q^n)$ is the isotropic
Lipschitz space on $Q^n$. As usual, the space
$\Lambda_\alpha(Q^n)$ is defined as the collection of all bounded
functions $g$ with compact support on $Q^n$ such that
\[
\sup_{(x,z),(y,w)\in Q^n}\Big\|g((x,z)\cdot
(y,w)^{-1})-g(x,z)\Big\|\le C\cdot \|(y,w)\|^{\alpha},\quad
0<\alpha<1.
\]
For $k<\alpha<k+1$ with $k=1,2,...$, we may define $g\in
\Lambda_\alpha(Q^n)$ is the set of all $C^k$ functions $g$ with
compact support on $Q^n$ such that $\nabla^k g\in
\Lambda_\alpha(Q^n)$.
\smallskip

\noindent $\bullet$ {\bf Sobolev spaces.} One may define the
non-isotropic Sobolev spaces $S^p_k(Q^n)$ with $k\in {\bf Z}_+$
and $1<p<\infty$ as follows
\[
S_k^p(Q^n)=\big\{f:Q^n\,\rightarrow \,\C:\,\, f\in L^p(Q^n), \,\,
{\mathcal P}(X,X^\prime)f\in L^p(Q^n)\big\},
\]
where ${\mathcal P}(X,X^\prime)$ is a monomial of degree $k$ in
vector fields $X,\,X^\prime\in V_1$. Here $L^p(Q^n)$ is the $L^p$
Lebesgue space.
\smallskip

\noindent $\bullet$ {\bf Hardy spaces.} The Hardy space $H^p(Q^n)$
with $0<p<\infty$ originally defined by maximal function as
follows: a distribution $f$ defined on $Q^n$ belongs $H^p(Q^n)$ if
and only if the maximal function
\[
{\mathcal M}(f)(x,z)=\sup_{\varepsilon>0}\big|f\ast
\phi_{\varepsilon}\big|(x,z)\,\in \, L^p(Q^n).
\]
Here $\phi\in {\mathcal S}(Q^n)$ with $\int_{Q^n}\phi(x,z)dxdz=1$
and polyradial (see Chapter 4 in Folland-Stein \cite{FS2}). As
usual,
\[
\phi_{\varepsilon}(x,z)=\varepsilon^{-(4n+6)}
\phi\big(\varepsilon^{-1}x,\varepsilon^{-2}z\big).
\]
The space $H^p(Q^n)$ can be defined by atomic decomposition and
maximal functions.

\begin{definition} A $H^p(Q^n)$ $p$-atom $(0<p\le 1)$ is a compactly
supported function $a(x,z)$ such that the following conditions
hold:

$(1)$ ({\it size condition}): there is a $Q^n$-ball
$B_{\bx_0}=B_{(x_0,z_0)}(r)=\{\bx=(x,z)\in Q^n:\,
\rho(\bx_0,\bx)\le r\}$ whose closure contains supp$(a)$ such that
$\|a(\bx)\|_{L^\infty}\le |B_{\bx_0}|^{-1/p}$;

$(2)$ ({\it moment condition}):
\[
\int_{Q^n}a(x,z){\mathcal P}(x,z)\,dxdz=0
\]
for all monomials ${\mathcal P}(\bx)={\mathcal P}(x,z)$ such that
${\mathcal P}(\bx)=\Big(\prod_{j=1}^n x_{j1}^{\alpha_{j1}}\cdots
x_{j4}^{\alpha_{j4}}\Big)z_1^{\beta_1}z_2^{\beta_2}z_3^{\beta_3}$
and
\[
\sum_{k=1}^4\sum_{j=1}^n\alpha_{jk}
+2(\beta_1+\beta_2+\beta_3)\le\Big[(4n+6)\big(\frac{1}{p}-1\big)\Big].
\]
Here $[s]$ is the integral part of $s$.
\end{definition}
Using the idea of atomic decomposition, we give the definition of
$H^p(Q^n)$ as follows
\[
H^p(Q^n)=\Big\{f\in {\mathcal S}^\prime(Q^n):\,
f=\sum_{k=1}^\infty \lambda_ka_k,\,\, {\mbox{where $a_k$ are
$p$-atoms}}, \,\,\sum_{k=1}^n |\lambda_k|^p<\infty\Big\},
\]
and we define
\[
\|f\|^p_{H^p(Q^n)}=\inf \Big\{\sum_{k=1}^n |\lambda_k|^p\Big\},
\]
where the infimum is taken over all possible atomic decompositions
of $f$. Then ``norm" $\|f\|^p_{H^p(Q^n)}$ is comparable to the
$\ell^p$ norm of the sequence $\{\lambda_k\}$ and the $L^p(Q^n)$
norm of the maximal function ${\mathcal M}(f)$.

Now, following notations in \cite{BGGR3}, for $g\in {\mathcal
S}(Q^n)$, one has
\begin{equation*}
\begin{split}
X_{jk}X_{j^\prime k^\prime}{\bf G}(g)(x,z)=&
\lim_{\varepsilon\rightarrow0^+}\int_{\rho(y,w)=\varepsilon}(X_{jk}G)(y,w)
(X_{j^\prime k^\prime}g)((x,z)\cdot(y,w)^{-1})\,dydw\\
& - \lim_{\varepsilon\rightarrow0^+}\int_{\rho(y,w)\ge
\varepsilon}(X_{jk}X_{j^\prime k^\prime}G)(y,w)
g((x,z)\cdot(y,w)^{-1})\,dydw.
\end{split}
\end{equation*}
The second term above is a generalized Calder\'on-Zygmund operator
in the sense of \cite{FS} and Koranyi-Vagi \cite{KV} because
$X_{jk}X_{j^\prime k^\prime}G$ is a kernel homogeneous of degree
$0$ satisfying mean value zero property. An operator ${\mathbf K}$
is said to be a generalized Calder\'on-Zygmund operator if the
following two conditions are satisfied:

\noindent {\rm(i)} ${\mathbf K}$ can be extended as a bounded
operator on ${L^2(Q^n)}$,

\noindent {\rm(ii)} There is a sequence of positive constant
numbers $\{C_j\}$ such that for each $j\in\mathbb N$,
\begin{equation}
\label{eq:ker1} \left (\int_{2^j\rho(\bx_2,\bx_3)\le
\rho(\bx_1,\bx_2)<
2^{j+1}\rho(\bx_2,\bx_3)}|K(\bx_1,\bx_2)-K(\bx_1,\bx_3)|^q\,d\bx_1\right
)^\frac{1}{q}\le C_j\cdot
\big(2^{j}\rho(\bx_2,\bx_3)\big)^{-\frac{4n+6}{q^\prime}}
\end{equation}
and
\begin{equation} \label{eq:ker2}
\left (\int_{2^j\rho(\bx_2,\bx_3)\le \rho(\bx_1,\bx_2)<
2^{j+1}\rho(\bx_2,\bx_3)}|K(\bx_2,\bx_1)-K(\bx_3,\bx_1)|^q\,d\bx_1\right
)^\frac{1}{q}\le C_j\cdot
\big(2^{j}\rho(\bx_2,\bx_3)\big)^{-\frac{4n+6}{q^\prime}},
\end{equation}
here $(q,q')$ is a fixed pair of positive numbers with
$\frac{1}{q}+\frac{1}{q^\prime}=1$ and $1<q^\prime<2$. We need the
following theorem to complete our discussion on the estimates for
${\bf G}$. The proof can be found in \cite{FS} and \cite{FS2}.
\begin{theorem} \label{th:2.1} Let ${\mathbf K}$ be a
generalized Calder\'on-Zygmund operator. Assume that the kernel
satisfies conditions (\ref{eq:ker1})-(\ref{eq:ker2}) with
$\{C_j\}\in \ell^1$. Then
\begin{equation*}
\|{\mathbf K}(f)\|_{L^p(Q^n)} \le C_p\|f\|_{L^p(Q^n)}, \quad
1<p<\infty;
\end{equation*} and
\begin{equation*} \Big|\{\bx\in Q^n:\, |{\mathbf
K}(f)(\bx)|>\lambda\}\Big|\le
\frac{C}{\lambda}\|f\|_{L^1(Q^n)},\quad \lambda>0.
\end{equation*}
\end{theorem}

When $0<p\le 1$, we may consider the boundedness of generalized
Calder\'on-Zygmund operator acting on Hardy spaces $H^p(Q^n)$. We
have the following theorem (see {\cite{C1}).
\begin{theorem}
\label{th:3.5} Let ${\mathbf K}$ be a generalized
Calder\'on-Zygmund operator. Suppose the following two conditions
hold:

\noindent {\rm(i)} (kernel assumption): there exist
$s\!\geq\big[(4n+6)(\frac 1{p}-\!1)\big]$ and $
\varepsilon\!>\!\frac1{p}\!-\!1$ such that
$$\{2^{(4n+6)j(2\gamma-1)}(C_j)^2\}\in \ell^1\quad\text{with}\quad
\gamma=1-\frac{1}{p}+\varepsilon,$$

\noindent {\rm(ii)} (adjoint operator assumption): ${\mathbf
K}^*(x_{1l}^{\alpha_{1l}}\cdots
x_{4l}^{\alpha_{4l}}z_1^{\beta_1}z_2^{\beta_2}z_3^{\beta_3})=0$
with $$\sum_{k=1}^4\sum_{l=1}^n\alpha_{kl}+\sum_{m=1}^3\beta_i\le
\big[(4n+6)(\frac{1}{p}-1)\big ],$$ where ${\mathbf K}^*$ is the
adjoint operator of ${\mathbf K}$.

\noindent Then ${\mathbf K}$ can be extended as a bounded operator
from $H^p(Q^n)$ into $H^p(Q^n)$ for $0<p\le 1$.
\end{theorem}

We first make the following observation. For any polynomial
${\mathcal P}(X)$ of degree $k$ in the horizontal vector fields
$X_{1l},\ldots,X_{4l}$ with $l=1,\ldots,n$, there exists another
polynomial $\tilde {\mathcal P}(X)$ such that the following
identity holds:
\begin{equation}
\label{eq:comm} {\mathcal P}(X){\bf G}={\mathbf G}\tilde {\mathcal
P}(X).
\end{equation}
Now we may apply the above theorems and (\ref{eq:comm}) to
conclude the following result.
\begin{theorem}
\label{th:est1} The fundamental solution $G(x,z)$ for the operator
$\Delta_0$ defines an operator ${\bf G}$ which satisfies the
following sharp estimates:

$i)$ $XX^{\prime}{\mathbf G}$ defines a bounded operator from the
Sobolev space $S^p_k(Q^n)$ into Sobolev space $S^p_k(Q^n)$ for
$k\in {\bf Z}_+$, $1<p<\infty$ and for all $X,\,X^{\prime}\in
V_1$;

$ii)$ $Z{\mathbf G}$ defines a bounded operator from the
non-isotropic Sobolev space $S^p_k(Q^n)$ into $S^p_k(Q^n)$ for
$k\in {\bf Z}_+$, $1<p<\infty$ and for all $Z\in V_2$;

$iii)$ $XX^{\prime}{\mathbf G}$ defines a bounded operator from
$H^p_k(Q^n)$ into $H^p_k(Q^n)$ for $k\in {\bf Z}_+$ and $0<p\le 1$
and for all $X,\,X^{\prime}\in V_1$;

$iv)$ $Z{\mathbf G}$ defines a bounded operator from $H^p_k(Q^n)$
into $H^p_k(Q^n)$ for $k\in {\bf Z}_+$ and $0<p\le 1$ and for all
$Z\in V_2$;

$v)$ ${\mathbf G}$ defines a bounded operator from
$\Lambda_\beta(Q^n)$ into $\Gamma_{\beta+2}(Q^n)\cap
\Lambda_{\beta+1}(Q^n)$ for $0<\beta<\infty$.

In conclusion, ${\mathbf G}$ gains two in horizontal directions
and only gains one in missing directions.
\end{theorem}

\end{document}